\documentclass[11pt]{amsart}
\def\D{{\mathbb D}}

\def\R{{\mathbb R}}

\newcommand{\ir}{\int}
\newcommand{\ep}{\varepsilon}
\newcommand{\eps}{\varepsilon}

\usepackage{color}
\usepackage{amsmath}
\usepackage{psfrag}
\usepackage{graphicx}
\usepackage{amssymb}
\usepackage{mathdots}
\usepackage{wasysym}
\usepackage{subfigure}
\usepackage{url}
\usepackage{cite}
\usepackage{array}
\usepackage{mathtools}

\usepackage{amsfonts,amssymb,amsmath,mathrsfs,graphicx,psfrag,color}
\usepackage{float}
\usepackage{epsfig}
\usepackage{graphicx}
\usepackage{caption}
\usepackage{kotex}
\usepackage{tikz}
\theoremstyle{plain}

\newtheorem{theorem}{Theorem}

\newtheorem{definition}[theorem]{Definition}

\newtheorem{remark}[theorem]{Remark}

\newtheorem{example}[theorem]{Example}

\numberwithin{equation}{section}
\numberwithin{theorem}{section}
\title[Heterogeneous nonlocal diffusion]{On the modelling of spatially heterogeneous nonlocal diffusion: deciding factors and preferential position of individuals}

\begin{document}

\maketitle

\begin{center}
{\large\bf Matthieu Alfaro \footnote{Universit\'e de Rouen Normandie, CNRS, Laboratoire de Math\'ematiques Rapha\"el Salem, Saint-Etienne-du-Rouvray, France \& BioSP, INRAE, 84914, Avignon, France.}, Thomas Giletti \footnote{IECL, Universit\'{e} de Lorraine, B.P. 70239, 54506
Vandoeuvre-l\`{e}s-Nancy Cedex, France.}, Yong-Jung Kim \footnote{Department of Mathematical Sciences, KAIST, 291 Daehak-ro, Yuseong-gu, Daejeon, 34141,
Korea.}, Gwena\"el Peltier \footnote{IMAG, Univ. Montpellier, CNRS, Montpellier, France.} and Hyowon Seo \footnote{Department of Applied Mathematics and the Institute of Natural Sciences, Kyung Hee University,
Yongin, 17104, South Korea.}}\\
[2ex]
\end{center}

\begin{abstract}
We develop general heterogeneous nonlocal diffusion models and investigate their connection to local diffusion models by taking a singular limit of focusing kernels. We reveal the link between the two groups of diffusion equations which include both spatial heterogeneity and anisotropy. In particular, we introduce the notion of deciding factors which single out a 
nonlocal diffusion model and typically consist of the \emph{total jump rate} and the \emph{average jump length}. In this framework, we also discuss the dependence of the profile of the steady state solutions on these deciding factors, thus shedding light on the preferential position of individuals.\\


\end{abstract}

%\tableofcontents

\section{Introduction}

The purpose of this exploratory paper is to investigate the modelling  of 
heterogeneous nonlocal diffusion and the connection with local diffusion. 
Starting from the most general nonlocal heterogeneous diffusion, we will see how to obtain a large class of local heterogeneous diffusion equations by taking a {\it focusing kernel limit}. This process sheds new light on the underlying dynamics of these models. In particular, the notion of {\it deciding factors}, and its implications on the shape of the steady state solutions, may help to find out whether a heterogeneous diffusion model is appropriate or not for a given application.

\subsection{Homogeneous diffusion} Under a homogeneous environment, a local diffusion equation is often given by the classical heat equation
\begin{equation}\label{heat}
u_t=D\Delta u,
\end{equation}
where the diffusivity $D>0$ is constant.  This diffusion equation is used 
in a wide range of applications from physics problems to ecology ones, where the homogeneous diffusion, $D \Delta u$, describes a random dispersal 
and $x\in\R^N$ the spatial position. Such a linear diffusion model can be 
recovered from a probabilistic individual based model. In the study of evolutionary biology, the linear diffusion also models the mutation process, where $x\in \R^N$ denotes a phenotypic trait (see \cite{Alf-Car-14, Fis-37,Kim-65,Kol-Pet-Pis-37,Shi-Kaw-97,Ske-51} among many other references). Remember that the diffusion equation \eqref{heat} is for an environment where there is no advantageous direction (isotropic) nor location (homogeneous). In ecology, evolutionary biology, epidemiology, etc. the model 
usually comes with a nonlinear reaction term, which takes into account birth, death, competition, etc. Throughout this paper, we neglect the nonlinear population dynamics and focus on the role of diffusion.

When the spatial domain is the whole Euclidean space $\R^N$, the nonlocal 
counterpart of the heat equation \eqref{heat} is
\begin{equation}\label{nonlocal_hom}
u_t =  \int_{\R^N} K( x-y) u(t,y) \, dy - m  u (t,x) ,
\end{equation}
where $ K : \R^N \to \R^+$ is a nonnegative integrable kernel and the coefficient $m$ is the total jump rate,
$$
m =\int_{\R^N} K(z) \, dz >0.
$$
The way to recover the local equation \eqref{heat} from the nonlocal one \eqref{nonlocal_hom} through a {\it focusing kernel limit} is classical and will be made precise in Section \ref{sec:singular}. The integro-differential equation \eqref{nonlocal_hom} is motivated by the long-range dispersal phenomenon meaning that an individual may ``travel'' long distance in a short time scale.  Nonlocal diffusions are therefore very popular in population models for which  long-distance dispersal events are relevant. 
The propagation of epidemics \cite{Mur-01,Med-Kot-03} and dispersal of seeds \cite{Cla-98} are such examples. Note that, in evolutionary biology, the mutation is also typically modeled in terms of integral operators \cite{Alf-Gab-Kav,Bur-86, Bur-88,Fle-79,Gil-17,Kim-65,Lan-75}. In the following, we mainly refer to \lq\lq spatial dispersal'' but, obviously, our analysis also applies to \lq\lq mutations''.

In the context of a spatial dispersal, the kernel $K (z)$ stands for the {\it jump rate} for an individual to move from $x$ to $x +z$ or, equivalently, $K(x-y)$ is the jump rate from $y$ to $x$. Hence, the constant $m$ is the {\it total jump rate} from a given point $x\in\R^N$. In particular, the last term in \eqref{nonlocal_hom} models the amount of species departing the position $x$. By integrating \eqref{nonlocal_hom}, one may check that the total population of the solution is preserved, i.e. $\int_{\R^N} u(t,x) dx$ is constant.

The central underlying assumption in the nonlocal diffusion model \eqref{nonlocal_hom} is  that the jump rate from $x$ to $x+z$ is independent from the departure point $x$ and the arrival one $x+z$. In most mathematical 
literature, the constant $m$ is normalized to 1 which can be done without 
loss of generality and the function $K$ is assumed to be radial (see \cite{Alf-Cov-17,Cha-Cha-Ros-06,Cov-Dup-07,And-10,Gar-11} and references therein).  As a result, the obtained nonlocal equation \eqref{nonlocal_hom} is isotropic and homogeneous in space, as is \eqref{heat}.

\subsection{Heterogeneous diffusion} However, the motility of biological species (in both local and nonlocal cases) depends on the environment, which is usually spatially heterogeneous. This spatial heterogeneity gives us a challenge to develop diffusion models that take into account the spatial heterogeneity properly in local and nonlocal equations. One of our goals is to address both issues jointly and to draw new connections between those models.

To give the motivation, let us first illustrate the challenge in modelling heterogeneous diffusion starting with a local diffusion case. A naive approach would be to simply replace the diffusivity $D$ in \eqref{heat} by 
a real-valued function $D(x)$ and obtain ``$u_t=D(x)\Delta u$". However, the resulting evolution equation does not satisfy the mass conservation 
property (even under the zero-flux boundary condition) and one should instead consider a diffusion model in a \emph{divergence form} for the motion of biological organisms.

On the other hand, there are infinitely many choices of diffusion equations in a divergence form when the diffusivity $D$ is not constant, which can be written 
\begin{equation}\label{qlaw}
u_t =\nabla\cdot\big(D(x)^q\nabla (D(x)^{1-q}u)\big),\quad q\in \R.
\end{equation}
Obviously, when the diffusivity coefficient $D$ is constant, these diffusion laws are all equivalent to the heat equation~\eqref{heat}. In a heterogeneous case, they lead to different solutions. Furthermore, if $D$ is not smooth enough, the \lq\lq classical weak solution'' is not defined and 
a different approach, depending on the choice of the diffusion law, is required.

Many well-known diffusion equations are written in the form of \eqref{qlaw}. First, one recovers the so-called Fick's diffusion law \cite{Fick1855} when $q=1$:
$$
u_t = \nabla \cdot \big ( D(x)  \nabla u)),\eqno{\text{(Fick's law)}}.
$$
This is a case without an advection phenomenon and constants are steady states. When $q=0$, the equation is written as
$$
u_t = \Delta ( D (x) u).  \eqno{(\text{Chapman's law})}
$$
Chapman derived this diffusion law using the kinetic theory in~\cite{Chapman1928}  where his main purpose is to explain the thermal diffusion phenomenon. Before Chapman's derivation, Wereide \cite{Wereide1914} derived a 
diffusion law that corresponds to the case with $q=\frac12$,
$$
u_t = \nabla \cdot \big ( \sqrt{D( x)} \nabla (\sqrt{D(x)} u ) \big),\eqno(\text{Wereide's law})
$$
using basic physical intuitive relations. The latter two equations are also respectively called It\^o and Stratonovich type since they are satisfied by probability density functions of the corresponding stochastic interpretations of the Brownian motion. Unlike Fick's law, the two contain drift phenomena and hence they may also be called Fokker-Planck equations.

The diffusion laws in \eqref{qlaw} satisfy the mass conservation property, and it is not obvious to see which one is the right diffusion law among all possibilities. The related question of whether heterogeneous dispersal is advantegeous to a species has been adressed from the evolutionary point of view in several works~\cite{MR3592434,MR2290095,CantrellCL}. Recently, two of the authors of the present work derived heterogeneous diffusion equations from revertible kinetic systems and compared the properties of diffusion laws using a thought experiment~\cite{KimSeo}. In particular, they obtained Wereide's diffusion law as a special case. In this paper we see that the same questions and issues exist and need to be addressed for the spatially heterogeneous nonlocal diffusion case. In this process we hope to obtain new insights on the long standing issue of heterogeneous diffusion (see, e.g., \cite{Nov-14} and the references therein).

\subsection{Organization of the paper} In Section~\ref{sec:nonlocal_het}, 
we propose a general model for nonlocal heterogeneous diffusion, introducing in particular the notion of {\it deciding factors}, from which a wide 
range of nonlocal heterogeneous diffusion laws can be derived. Next, in Section~\ref{sec:singular}, we draw a connection between the nonlocal and the local diffusion cases by a singular limit argument. This provides some new interpretation for equation~\eqref{qlaw}. However, as was pointed out by two of the present authors in an earlier work~\cite{KimSeo}, the diffusivity function~$D(x)$ may not be sufficient to describe heterogeneous 
diffusion accurately, because it does not fully cover velocity jump processes given by kinetic equations. In subsection~\ref{subsec:double}, we make a related observation that the nonlocal diffusion equation, or position jump process, also leads to a range of local equations wider than \eqref{qlaw}. In Section \ref{sec:steady}, in order to answer the question \lq\lq Where are the individuals?'', we investigate steady states solutions (through both rigorous computations and numerical simulations) of the equations with heterogeneous nonlocal diffusion. This in particular reveals the role of deciding factors on the repartition of individuals, and may guide the use of these models in the applications. Finally, in Section \ref{s:summarize}, we briefly summarize our main contributions and present some perspectives.

\section{Heterogeneous nonlocal diffusion equations}\label{sec:nonlocal_het}

In a nonlocal ecological diffusion model, it is assumed that individuals have the ability to jump from one position to another at any given time. In particular, we occasionally refer to the corresponding equation as the 
(position) jump process in order to distinguish it from kinetic equations 
which take integro-differential terms of velocity jumps. As mentioned earlier, a commonly used jump process is the homogeneous case \eqref{nonlocal_hom}, where the jump rate depends only on the distance between two points.

However, if the environment is heterogeneous, the jump rate may depend on 
the departing point or arrival one. We may denote the jump rate by
$$
J (x,y): \text{ the  jump rate from $x$ to $y$}.
$$
Let us emphasize again that this jump rate depends not only on the distance between the two points, but also on the direction of the jump. Thus the spatial heterogeneity of the environment along the way is counted by considering the distance {\it and} direction. This leads to a general nonlocal heterogeneous equation
\begin{equation}\label{nonlocal_het_gen}
u_t (t,x) = \int_{\R^N} J(y,x) u(t,y) \, dy - \int_{\R^N} J (x,y)  u(t,x) \, dy.
\end{equation}
Here, the first integral term stands for individuals arriving at $x$ from 
all possible positions, whereas the second one accounts for individuals leaving from the position $x$. We refer to \cite{Hut-Mar-Mis-Vic-03} for its derivation from a discrete model in space and time. Notice that, provided that $J$ is integrable with respect to each variable, the above equation is well-posed and satisfies the mass conservation property (at least formally) by Fubini's theorem. 

Observe that~\eqref{nonlocal_het_gen} is recast
$$
u_t(t,x)=\int_{\R^N} J(y,x) u(t,y) \, dy - m(x)u(t,x),
$$
where 
\begin{equation}
\label{def:mx}
m(x):=\int_{\R^{N}} J(x,y)dy
\end{equation}
stands for the total jump rate from position $x$. For late use, we also define
\begin{equation}
\label{def:J1}
J_1(x):=\int_{\R^N} \Vert y-x\Vert J(x,y)dy
\end{equation}
the absolute first moment at position $x$. In other words, $J_1(x)$ stands for the average jump length among individuals jumping from the position 
$x$. Both the total jump rate and the average jump length are natural quantities which may be impacted by spatial heterogeneity.

Obviously, equation \eqref{nonlocal_het_gen} turns into the  homogeneous nonlocal diffusion equation \eqref{nonlocal_hom} by selecting $J(y,x) =K(x-y)$ and defining $m := \int_{\R^N} K (z) dz$. Yet our purpose is now to find practical ways to take into account the spatial heterogeneity and anisotropy in the choice of the jump rate function~$J$. We will see that it is advantageous to specify~$J$ by distinguishing the effects of heterogeneity on  the distance and on the direction of jumps.

\subsection{A few illuminating examples} Before we introduce general ways 
to handle heterogeneity, we start with some revealing examples. First, in 
the context of  evolutionary biology, let us mention the work \cite{Gil-19} which considers a mathematical model for the fitness distribution in an asexual population under mutation and selection. Due to the presence of 
a phenotype optimum, the distribution of mutation effects on fitness depends on the parent’s fitness (corresponding to $x$). Hence the authors need to consider  \lq\lq context-dependent" mutation kernels, corresponding in our setting to the case where 
\begin{equation}\label{eq:departure}J(x,y) = K(x;y-x).
\end{equation}
In other words a dispersal kernel $K(x;\cdot)$ has been assigned to each point~$x$ in the spatial domain, accounting for heteregeneity. In this above form, $J(x,y)$ is the jump rate from $x$ to $y$, and therefore the heterogeneity factor is considered at the {\it departure point}. This case arises from an It\^o interpretation of a stochastic Poisson jump process, 
and therefore by analogy with the local diffusion we may call it nonlocal 
It\^o type diffusion.

To make the above form more meaningful, one may further specify the kernel function~$K$, for instance by letting
\begin{equation}\label{eq:departure_m}
K(x;y-x)= m(x) \widetilde{K} (y-x),
\end{equation}
or 
\begin{equation}\label{eq:departure_g}
K(x;y-x) =\frac{1}{g(x)^N} \widetilde{K} \left( \frac{y-x}{g(x)} \right),
\end{equation}
where in either cases $\widetilde{K}$ denotes a probability density, which we may rescale without loss of generality so that its absolute first moment
$$\int_{\R^N} \| z\| \widetilde{K} (z) dz = 1.$$
In the former case, the $x$-dependence appears only in the factor $m(x)$ which as outlined before stands for the total jump rate from the position~$x$. In the latter case, one may compute that the total jump rate from~$x$ is constant equal to~$1$, while the average jump length (or absolute first moment at position~$x$) is $g(x)$ and now accounts for the heterogeneity. Both of these forms provide a tractable way to include heterogeneity in the model.

However, depending on the context, one may want to consider a symmetrical 
situation where the jump rate from~$x$ to~$y$ is the same as the jump rate from~$x$ to~$y$. The reasoning behind that may be that the jump rate of 
individuals may only depend on the path of the displacement, and not necessarily on heterogeneous conditions at the departure point. The above examples are ill-suited for such a situation, and leads us to instead consider the form
\begin{equation}\label{eq:middle}
J(x,y) = K \left( \frac{x+y}{2}; y-x \right).
\end{equation}
The symmetry of $J(x,y)$ can then be reformulated as the evenness of the kernel $K (p; \cdot)$, for any $p \in \R^N$. Here $p = \frac{x+y}{2}$ stands for the point where heterogeneity dependence is involved: in this framework, heterogeneity is taken at the {\it middle point of the jump}.

The common trait to the above two forms~\eqref{eq:departure} and~\eqref{eq:middle} is that we  assign a kernel to each position in the space domain. Then the heterogeneity corresponding to the path $(x,y)$ is determined 
by environmental conditions at a single point~$p$, which is then called a 
\textit{deciding factor}. Extrapolating on the above examples, we are led 
to pick
$$
p = \alpha x + \beta y,
$$
with $\alpha + \beta = 1$, in other words as a barycenter of the departure and arrival points. When $\alpha =1$, we find the nonlocal It\^o type diffusion where the influential environmental conditions are those at the starting point~$x$. Let us  notice that this particular case  has been considered in~\cite{Cor-Cov-Elg-07,Cor-Elg-Gar-Mar-09,Cor-Elg-Gar-15, Cor-Elg-Gar-16}. When $\alpha = \frac{1}{2}$ we recover the symmetrical case when the jump rates are the same from both $x$ to $y$ and $y$ to $x$. A key input of this work is to allow any barycenter of the starting and arrival points to be a deciding factor. This will be discussed again in subsection~\ref{ss:deciding}.

Still, it turns out that a single deciding factor may not be enough to fully understand the diffusion process, and we provide one last example to illustrate this. Let us consider, in spatial dimension $N=1$, the case
\begin{equation}\label{eq:strat}
J (x,y ) = \frac{1}{h(y)} \widetilde{K} \left( \int_x^y  \frac{1}{h (s) 
} ds \right),
\end{equation}
where, as above, $\widetilde{K}$ is a probability density and $h$ is positive. This nonlocal model arises from a Stratonovich (also referred to as 
Marcus~\cite{applebaum_2009}) integral of a Poisson jump process~\cite{Bartlett}. Another way to derive it is to check that the heterogenous  equation~\eqref{nonlocal_het_gen} with~\eqref{eq:strat} reduces to the homogenenous nonlocal diffusion model~\eqref{nonlocal_hom} by to the change of variables
$$x' = \int_0^x \frac{1}{h(s)} ds.$$
This is consistent with the property that the Stratonovich integral preserves the chain rule of standard calculus. From the biological point of view, this means that the heterogeneity distorts the distances, as in the jump rate the (Euclidean) distance between $x$ and $y$ is replaced by the integral from~$x$ to~$y$ of the inverse of the function~$h$. This integral may be understood as a distance with respect to a so-called {\it food metric}~\cite{ChoiYJK}.

Due to the integral term in \eqref{eq:strat}, one cannot isolate a single 
deciding point and it is the {\it whole path} which determines the heterogeneity. This partly explains why, in the next subsection, we will introduce a general framework where a kernel is assigned not only to any point (though, as suggested by the above examples, it is sufficient in some situations) but to any possible path. Such generality also allows us to consider the case of two deciding factors, which is convenient to distinguish 
the effects of heterogeneity on  the total jump rate and on the average jump length. 

\subsection{Assignment of a kernel to any possible path}  In order to encompass all the above examples in a very general setting, our idea is the following: for any possible path, i.e. for any pair of points $(x,y)$, we 
{\it assign} a dispersal kernel $K(x,y;\cdot)$. Notice that the case $K(x,y;\cdot)=K(x,\cdot)$ corresponds to the situation in~\cite{Gil-19} mentioned above.

Precisely, we assume the form
$$
J (x,y) =K (x, y ; y-x),
$$
where
$$ K : (x,y ;z) \in \R^N \times \R^N \times \R^N \mapsto K(x,y;z)\in \R,$$
is a nonnegative function such that, for any $(x,y) \in \R^N \times \R^N$,
$$
m(x,y) := \int_{\R^N} K (x,y;z) \, dz  \in (0, +\infty),
$$
and
$$
J_1(x,y):=\int_{\R^{N}}\Vert z\Vert K(x,y;z)\, dz  \in (0, +\infty).
$$
Notice that, in particular, to any pair of points $(x,y)$ there \lq\lq corresponds'' a nonlocal equation of type~\eqref{nonlocal_hom}, where $m (x,y)$ and $K ( x ,y ; \cdot)$ are respectively the total jump rate and the 
dispersal kernel. Such a \lq\lq fictive'' equation could be interpreted as the diffusion law in a homogeneous medium where environmental conditions are everywhere the same as they are perceived by individuals jumping from $x$ to $y$. In the sequel, even if \lq\lq possibly fictive'', we denote $m(x,y)$ the total jump rate and $J_1(x,y)$ the absolute first moment, corresponding to the path $(x,y)$.

With the above choice of the jump rate $J$, the general equation~\eqref{nonlocal_het_gen} is recast
\begin{equation}\label{nonlocal_het}
u_t (t,x) = \int_{\R^N} K (y, x; x-y )  u(t,y) \, dy -  \int_{\R^N}  K ( x,y ; y -x) u (t,x) \, dy.
\end{equation}
To ensure that the right hand side of the above equation is  \lq\lq meaningful'',  we will also always assume that there exists $\overline{K} \in L^1 (\R^N)$ such that
$$
\forall (x,y) \in \R^N \times \R^N, \quad K ( x,y ; \cdot) \leq \overline{K}.
$$

Equations~\eqref{nonlocal_het_gen} and~\eqref{nonlocal_het}, together with the above assumptions, provide a general heterogeneous framework within 
which all the following sections will fit. Obviously, the latter equation~\eqref{nonlocal_het} still includes the homogeneous equation: it naturally arises when we associate the same dispersal kernel to all paths, i.e. regardless of the departure and arrival points. It also includes \eqref{eq:strat}, by letting
$$
K(x,y;z) = \frac{1}{h(y)} \widetilde{K} \left( z  \int_0^1 \frac{1}{h(x+ s(y-x))} ds \right),
$$
provided that $0 < \inf_{\R^N} h \leq \sup_{\R^N} h < +\infty$.\medskip

In order to emphasize the role of the total jump rate and the absolute first moment, we point out that a natural generalization of example~\eqref{eq:departure_m} consists in
\begin{equation}\label{ex1}
K (x,y ;z ) = m(x,y) \widetilde{K} (z),
\end{equation}
where $\widetilde{K}$ is a probability density. In this case, the heterogeneity only appears in the total jump rate $m(x,y)$ but not in the actual 
jump distribution $\widetilde{K}$.

Similarly, a natural  generalization of example~\eqref{eq:departure_g} consists in
\begin{equation}\label{ex2}
K (x, y ;z ) = \frac{1}{g(x,y)^N} \widetilde{K} \left( \frac{z}{g(x,y) } \right),
\end{equation}
where $\widetilde{K}$ is a  probability density, rescaled so that its absolute first moment
\begin{equation}\label{absolute}
\int_{\R^N} \|z \|  \widetilde{K} (z) \, dz = 1,
\end{equation}
and the function $g : \R^N \times \R^N \to \R$ has positive infimum and supremum.
Moreover, the total jump rate
$$
m(x,y) = \int_{\R^N}  K (x,y ;z ) \, dz=\int_{\R^N}\widetilde K (z)\,dz
 $$
is actually spatially constant. However, the heterogeneity appears in the 
function $g(x,y)$, which is the (fictive) {\it relative absolute first moment} of the dispersal kernel $K (x,y ; \cdot)$, see \eqref{relative} below, and can be interpreted as a (fictive) {\it average jump length}. 

Finally, we may combine \eqref{ex1} and \eqref{ex2} so that both the total jump rate and the average jump length are spatially heterogeneous. The heterogeneous dispersal kernel can then be written as
\begin{equation}\label{typical}
K ( x,y ;z) = \frac{m (x,y)}{g(x,y)^N} \widetilde{K} \left( \frac{z}{g(x,y) } \right),
\end{equation}
where $\widetilde{K}$ is a  probability density, rescaled so that its absolute first moment is equal  to one. The relative absolute first moment of the dispersal kernel $K(x,y;\cdot)$ is defined by
\begin{equation}
\label{relative}
\frac{\displaystyle{\int_{\R^N} \|z \|  K(x,y;z)\, dz}}{\displaystyle{\int _{\R^{N}} K(x,y;z)\,dz}}
\end{equation}
and is $g(x,y)$ when assuming form \eqref{typical}.

\subsection{The notion of deciding factors}\label{ss:deciding} In Section 
\ref{sec:singular}, we will compute the singular limits of \eqref{nonlocal_het} with focusing kernels to recover a wide range of local diffusion equations, including, but not limited to, \eqref{qlaw}. The point of such focusing kernels is that the  average jump length goes to 0 as a singular 
limit parameter $\varepsilon>0$ tends to 0. In this context, it turns out 
that it is mainly sufficient that the distribution kernel depends on a finite number of points. For instance, in the above example \eqref{typical}, one may assume that
$$g (x,y) = \widetilde{g} (\alpha x + \beta y),$$
$$m (x,y) = \widetilde{m} (\alpha ' x + \beta ' y),$$
where $\alpha + \beta = \alpha ' + \beta ' = 1$, so that
\begin{equation*}
J(x,y) = \frac{\widetilde m (\alpha'x+\beta'y)}{(\widetilde g(\alpha x+\beta y))^N} \widetilde{K} \left( \frac{y-x}{\widetilde g(\alpha x+\beta y) } \right).
\end{equation*}
 This means that the total jump rate~$m$ of the dispersal kernel corresponding to the path $(x,y)$ is determined by environmental conditions at the single point $\alpha' x + \beta' y$, which as outlined before we then call a {\it deciding factor}. On the other hand, the relative absolute first moment of the dispersal kernel corresponding to the path $(x,y)$, which is also the average jump length $g(x,y)$, is determined by environmental conditions at the (possibly) distinct deciding factor $\alpha  x + \beta y$. It is left to decide how to choose those two points, and to observe 
the consequences in the singular limit depending on that choice. For instance, one may pick $\alpha =1 $ and $\beta =0$, so that the average jump length $g(x,y)$ is decided at the departure point $x$, while $\alpha ' = \beta ' = \frac{1}{2}$, which means that the total jump rate $m(x,y)$ is decided by the middle point $\frac{x+y}{2}$. We believe that this 
interpretation, together with the connection between local and nonlocal diffusion which we next establish, provide a new perspective on the understanding of heterogeneous diffusion.

\section{Singular limits: from nonlocal to local diffusion}\label{sec:singular}

In this section, we show how to obtain a large class of local diffusion laws from heterogeneous jump processes of the type \eqref{nonlocal_het}. Recall that, in the homogeneous case, one can recover the heat equation \eqref{heat} as the {\it focusing kernel limit} of the integro-differential 
\eqref{nonlocal_hom}. Precisely, assuming that $K$ has a finite second moment the following holds: as $\ep \to 0$, the solution $u_\ep$ of \eqref{nonlocal_hom}, with the kernel  $K(x)$ replaced by the focused kernel
$$
K_\ep(x):=\frac{1}{\ep^{N+2}}K\left(\frac x \ep\right),
$$
tends (in a sense we do not precise in this exploratory paper) to the solution of \eqref{heat} with $D$ depending on the second moment of $K$, starting with the same initial datum. This can be understood from a simple formal Taylor expansion that we will indicate below in a much more general 
setting. Hence, in some situations, nonlocal dipersal operators can be approximated by the local diffusion operator. This fact was long known: for 
instance, in evolutionary biology models, we refer to \cite{Kim-65} and \cite[Chapter VI, subsection 6.4]{Bur-00-book}.  For a rigorous proof, one 
can use the explicit writing of the solution of \eqref{nonlocal_hom} in the Fourier variable, see \cite[Theorem 1.24]{And-10}.

As far as the singular limit of heterogeneous nonlocal diffusion equations is concerned, we refer to the recent works \cite{Mol-Ros-16, Mol-Ros-19,Sun-Li-Yan-11} for some particular cases that all fall into our general framework below. Let us also notice the preprint~\cite{Dos-Oli-Ros-20} which starts from a local/nonlocal heterogeneous model.

In view of the very general framework outlined in Section~\ref{sec:nonlocal_het}, hereafter we will use the focusing kernels
$$K_\ep (x,y; z) := \frac{1}{\ep^{N+2}} K\left(x,y; \frac{z}{\ep} \right).$$
Thanks to the \lq\lq additional'' variable $z$ (ultimately replaced by $y-x$ in the diffusion equation), the above form drives the average jump length to 0 as $\ep \to 0$, while preserving the scale of the heterogeneity 
in the singular limit.

\subsection{The diffusivity matrix}The limiting local equation will involve the following notion of diffusivity.

\begin{definition}[The diffusivity matrix]\label{diffusivity}
Let a dispersal kernel $K (x,y; z)$ as in  Section~\ref{sec:nonlocal_het}. For all $p \in \R^{N}$, we define a diffusivity matrix by
\begin{equation}\label{def:D}
\D(p):=\frac {1}{2} \ir  _{\R^N} z\otimes z  K (p,p;z) \, dz.
\end{equation}
Alternatively, the diffusivity matrix $\D(p)=(d_{ij}(p))_{1\leq i,j \leq n}$ can be written component-wise as
\begin{equation}\label{def:dij}
d_{ij}(p)=\frac {1}{2}\ir _{\R^N}   z_iz_j K(p,p;z) \, dz.
\end{equation}
\end{definition}

Notice that in the integral terms of the above definition, we only need to evaluate $K(p,p';z)$ at $p=p'$. The reason is that the local equation 
arises in the limit of focusing kernels where the average jump length goes to 0, hence the departure and arrival point eventually coincide. There is however no mathematical obstacle to extend the above diffusivity matrix to the case where $p \neq p'$. Notice also that it is implicitly assumed that the kernel $K$ is such that  $d_{ij}(p)\in \R$ for any~$p \in \R^N$.

In order to illustrate this notion, let us consider the aforementioned example \eqref{typical}, namely
$$
K ( x,y ;z) = \frac{m (x,y)}{g(x,y)^N} \widetilde{K} \left( \frac{z}{g(x,y) } \right),
$$
for which the diffusivity matrix is straightforwardly computed as
$$
\D (p) := m(p,p) g(p,p)^2 \widetilde{\D}.
$$
Here, $\widetilde{\D}$ is the constant diffusivity matrix associated with 
the homogeneous dispersal kernel $\widetilde{K}$, while the total jump rate $m$ and the average jump length~$g$ account for the environment's heterogeneities. Unsurprisingly, the diffusivity depends monotonically on both the total jump rate and the average jump length.

\subsection{The case of a single deciding factor}\label{subsec:one-point} 
As outlined at the end of Section~\ref{sec:nonlocal_het}, we consider the 
nonlocal heterogeneous model \eqref{nonlocal_het} with
$$K ( x ,y; \cdot) =  \widetilde{ K} (\alpha x + \beta y ; \cdot),$$
where
$$\alpha + \beta =1.$$
In other words, the nonlocal dispersal kernel associated with a path $(x,y)$ is determined by the single point $\alpha x + \beta y$. We immediately drop the tilde for convenience and rewrite~\eqref{nonlocal_het} as
\begin{equation}\label{nonlocal_het_1}
u_t (t,x) = \int_{\R^N} K (\alpha y + \beta x ; x-y)  u(t,y) \, dy -  u 
(t,x) \int_{\R^N} K (\alpha x + \beta y ; y -x)  \, dy.
\end{equation}
We point out that, when $\alpha =0$ and $\beta=1$, the kernel is chosen from the arrival point and the equation is recast
\[
u_t(t,x)=  \ir_{\R^N}  K (x; x-y ) u(t,y)\,dy-u(t,x)\ir_{\R^N}   K(y; y-x)\,dy.
\]
When $\alpha =1$ and $\beta=0$, the kernel is chosen from the departing point and the equation is  recast
\[
u_t(t,x)=\ir_{\R^N}   K(y ; x-y) u(t,y)\, dy-  u (t,x) m(x) ,
\]
where as before $m$ denotes the integral of $K$ with respect to its last variable. When $\alpha\notin\{0,1\}$ the dispersal kernel is selected by a nontrivial linear combination of the departure and the arrival points. If we integrate \eqref{nonlocal_het_1} over $x\in \R^{N}$ and use Fubini's theorem, we obtain as explained earlier the mass conservation property:
$$
\frac{d}{dt}\ir_{\R^N} u(t,x)\,dx=0.
$$
Here we implicitly assumed the framework introduced in Section~\ref{sec:nonlocal_het}: the function $K (p ; \cdot)$ is nonnegative and bounded by some $\overline{K} \in L^1 (\R^N)$ uniformly with respect to $p \in \R^N$.

Moreover, we make the additional assumption that $K$ is symmetric with respect to its last variable, i.e.
\begin{equation}\label{hyp-sym}
\forall (p,z) \in \R^N  \times \R^N, \qquad K (p; -z) = K(p;z).
\end{equation}
We also assume that $K$ is \lq\lq sufficiently smooth'' in the $p$-variable, and that its derivatives with respect to the $p$-variable are \lq\lq sufficiently $z$-integrable'' (typically quadratic weights have to be supported). These assumptions play an important role in the singular limit argument, as will be clear from the terms appearing in the computations below.

Let us point out that a typical example of such dispersal kernels is
$$K (p;z) = m(p) \widetilde{K} (z),$$
where $\widetilde{K}$ is a symmetric probability density, and thus only the mass $m(p)$ of the dispersal kernel depends on position~$p$. One may also assume that the relative absolute first moment, or average jump length, is spatially heterogeneous by considering
$$K (p; z) = \frac{m(p)}{g(p)^N} \widetilde{K} \left( \frac{z}{g(p)} \right).$$
\medskip

We now formally derive local diffusion equations from the nonlocal diffusion equations \eqref{nonlocal_het_1}. To do so, we actually take a sequence of focusing kernels. Denote, for any $z \in \R^N$,
\begin{equation}\label{Keps}
K_\ep ( p; z):=\frac{1}{\ep^{N+2 }}K\Big(p; \frac z \ep  \Big),\quad 0<\ep \ll 1.
\end{equation}
Notice that these still satisfy our assumptions. However, with obvious notations, the total jump rate and the average jump length are related through
$$
m_\ep(p):=\int_{\R^N} K_\ep(p;z)\, dz=\frac{1}{\ep^2}m(p), \quad g_\ep(p):=\frac{\displaystyle{\int_{\R^N} \| z \|  K_\ep(p;z)\, dz}}{\displaystyle{\int _{\R^{N}} }K_\ep(p;z)\, dz}=\ep g(p).
$$
Hence, as $\ep \to 0$, the average jump length goes to 0 but, in some sense, this is compensated by the fact that $m_\ep(p) \to +\infty$, i.e. individuals jump more often. In particular, the corresponding diffusivity matrix $\D_\ep$ of the focusing kernels is unchanged. This can be understood from the fact that, in the homogeneous case, the kernels~\eqref{Keps} naturally arise from the change of variables
$$
t \leftarrow \varepsilon^2 t \, , \quad x \leftarrow \varepsilon x.
$$
Thus, in our context, focusing kernels can be understood as having the motion of individuals occur on a smaller scale than the heterogeneity.

Let us proceed with our computation of the singular limit $\ep \to 0$. First, we plug the above focusing kernels into \eqref{nonlocal_het_1}, use a change of variable $z=\frac{y-x}{\ep}$, and obtain
\begin{equation*}
\eps^2u_t(t,x)=\ir _{\R^N} K\Big( x+\ep\alpha z ;z \Big)u(t,x+\ep z)dz-\ir_{\R^N} K\Big( x+\ep\beta z ;z\Big)u(t,x)dz,
\end{equation*}
where we used the symmetry of the dispersal kernel with respect to its last variable. Then
\begin{eqnarray*}
&&\eps^2u_t(t,x) \approx \ir_{\R^N}  \Big(K(x;z)+\ep \alpha z\cdot \nabla_p K(x;z)+\ep^{2}\frac 1 2 \alpha^{2}\left\langle D^{2}_p K (x;z)z,z\right\rangle \Big)\\
&&\;\times\Big(u(t,x)+\ep z\cdot \nabla u(t,x)+\ep^{2}\frac 12 \left\langle D^{2 }u(t,x)z,z\right\rangle\Big)dz\\
&&\;-\ir _{\R^N}\Big(K(x;z)+\ep \beta  z\cdot \nabla_ p K(x;z)+\ep^{2}\frac 1 2 \beta ^{2}\left\langle D^{2}_p K(x;z)z,z\right\rangle \Big)u(t,x)dz,
\end{eqnarray*}
where $\nabla _p K (x;z)$ and $D^{2}_p K(x;z)$ respectively denote the gradient vector and the Hessian matrix of the map of dispersal kernels, $p\mapsto K(p;z)$.

Clearly the $\ep^{0}$ order terms in the right hand side cancel each other. Moreover, since $K(p;z)$ is symmetric with respect to $z$, we have
$$
\ir _{\R^N} z\cdot \nabla_p K(x;z)\,dz=0, \quad \ir_{\R^N} K(x;z)z\cdot 
v \,dz=0,
$$
for any $v\in \R^N$. Therefore, the $\ep^{1}$ order terms in the right hand side  also cancel out. As a result, recalling that $\beta=1-\alpha$, 
we reach the local equation
\begin{eqnarray*}
u_t(t,x)&\approx &\frac 12 \ir_{\R^N}  K(x;z) \langle D^2u(t,x)z,z\rangle\,dz\\
&&+\alpha\ir _{\R^N}\Big(z\cdot \nabla_p K(x;z)\Big)\Big(z\cdot \nabla u (t,x)\Big)\,dz\\
&&+\frac{2\alpha-1}{2}\ir _{\R^N}\langle D^2_pK(x;z)z,z\rangle u(t,x)\,dz.
\end{eqnarray*}

We now define the \lq\lq exponent''
\begin{equation}
\label{def:p}
q:= 2 - 2 \alpha \, (= 2\beta).
\end{equation}
Then, using the exponent $q$ and the coefficients $d_{ij}(p)$ of the diffusivity matrix from \eqref{def:dij}, we may rewrite the above diffusion equation as
\begin{eqnarray*}
u_t(t,x)&\approx & \sum_{i,j}d_{ij}(x)\partial_{x_ix_j}u(t,x)+(2-q)\sum_{i,j}\partial _{p_j}d_{ij}(x)\partial _{x_i}u(t,x)\\
&&\quad +(1-q)\Big(\sum_{i,j} \partial_{p_i p_j}d_{ij}(x)\Big)u(t,x).
\end{eqnarray*}
We pursue and, omitting to write the $(t,x)$ variable, get
\begin{eqnarray*}
u_t&\approx & \sum_{i,j}\partial_{x_j}\Big(d_{ij}(x)\partial_{x_i}u\Big) \\&&+(1-q)\Big(\sum_{i,j}\partial _{p_j}d_{ij}(x)\partial _{x_i}u+\Big(\sum_{i,j} \partial_{p_i p_j}d_{ij}(x)\Big)u\Big)\\
&\approx & \sum_{i,j}\partial_{x_j}\Big(d_{ij}(x)\partial_{x_i}u\Big) +(1-q)\sum_{i,j}\partial_{x_i}\Big(\partial _{p_j}d_{ij}(x)u\Big)\\
&\approx & \sum_{i}\partial_{x_i}\Big(\sum _j (d_{ij}(x)\partial_{x_j} u+(1-q)\partial_{p_j}d_{ij}(x)u)\Big),
\end{eqnarray*}
which is in the so-called divergence form. We finally obtain, in the singular limit, the  (possibly anisotropic) diffusion equation
\begin{equation}\label{componentwise}
u_t=\sum_{i}\partial_{x_i}\Big(\sum _j\Big(d_{ij}^{q}(x)\partial_{x_j}\Big(d_{ij}^{1-q}(x)u\Big)\Big)\Big),
\end{equation}
which is an anisotropic generalization of the form~\eqref{qlaw}.

\begin{remark}[Isotropic diffusion]
The relation between \eqref{componentwise} and~\eqref{qlaw} becomes obvious in the isotropic case, i.e. when radial symmetry is assumed for the mapping $z\mapsto K(p;z)$, which is a stronger assumption than~\eqref{hyp-sym}. In this case, we obtain $\D (p)= D (p)\mathbf{I}$, where $\mathbf{I}$ is the identity matrix and
$$
D(p):=\frac{1}{2N} \ir _{\R^N} \| z \|^2 K(p;z)\,dz.
$$
The diffusion equation~\eqref{componentwise} is then recast
\[
u_t=\nabla\cdot\Big(D^{q}(x) \nabla \Big(D^{1-q}(x)u\Big)\Big),
\]
which is exactly~\eqref{qlaw}.
\end{remark}

Let us first observe that, if $q=1$ in \eqref{componentwise}, we reach Fick's law:
\[
  u_t = \nabla\cdot(\D(x)\nabla u).\ \quad\eqno(q=1, \alpha =\frac12)
\]
This is the case when $\alpha=\beta=\frac12$, meaning that  the middle point between departure and arrival points is used to define the spatial heterogeneity in the nonlocal model (for instance, in a heterogeneous total jump rate $m$).

The next case is when $q=0$ in \eqref{componentwise}, which yields Chapman's law:
\[
  u_t = \nabla\cdot(\nabla\cdot(\D (x) u)). \eqno(q=0,\alpha=1)
\]
Here, $\nabla\cdot(\D u)$ is the vector that consists of the divergence of the rows of the matrix $\D u$, that is
$$
\nabla\cdot(\D u)=\Big(\nabla\cdot(d_{i1}u,...,d_{iN}u) \Big)_{1\leq i \leq N}.
$$
This is the case when $\alpha=1$, $\beta=0$, meaning that  only the departure point is used to define the spatial heterogeneity in the nonlocal model.

For a general $q$, we may rewrite  equation \eqref{componentwise} as follows. First, observe that
\[
\sum _j\partial _{x_j}\Big(d_{ij}\,u\Big)=
\sum _j d_{ij}^q\partial _{x_j}\Big(d_{ij}^{1-q}u\Big)+\sum_j \Big(\partial _{x_j}(d_{ij}^q) d_{ij}^{1-q}\Big)u.
\]
Let $N_q(x)=((N_q)_i(x))$ be the vector given by
\begin{equation}\label{Nq}
(N_q)_i(x)=\Big(\sum_j \partial _{x_j}(d_{ij}^q(x)) d_{ij}^{1-q}(x)\Big).
\end{equation}
Then, the component-wise equation \eqref{componentwise} can be written as
\begin{equation}\label{ANq}
u_t=\nabla\cdot\Big(\nabla\cdot(\D (x) u)-N_q(x)u\Big).
\end{equation}

\begin{remark}[Orthotropic diffusion]
Assume the component-wise symmetry, which is a stronger than the symmetry 
assumption \eqref{hyp-sym}: for all $(p,z)\in \R^{N}\times \R^{N}$,
\begin{equation}\label{sym-stonger}
 K(p; \alpha_1 z_1,...,\alpha_Nz_N)=K(p; z)\ \ \text{for all } (\alpha_i)\in \{\pm 1\}^N.
\end{equation}
Then, the diffusivity matrix $\D(p)$ becomes diagonal, which is often called orthotropic. Under assumption \eqref{sym-stonger}, equation \eqref{componentwise} turns into
\[
u_t=\sum_{i}\partial_{x_i}\Big(d_{ii}^{q}(x)\partial_{x_i} \Big(d_{ii}^{1-q}(x)u\Big)\Big).
\]
Furthermore, since the diffusivity matrix $\D(p)$ is diagonal, $\D^q(p)$ can be defined element-wise, and the equation is recast
\[
u_t=\nabla\cdot\Big(\D^{q}(x) \Big(\nabla \cdot (\D^{1-q}(x)u)\Big)\Big).
\]
\end{remark}

\subsection{The case of two deciding factors}\label{subsec:double} We now 
consider the nonlocal heterogeneous model \eqref{nonlocal_het} with
\begin{equation}\label{ass:separation}
K ( x ,y; \cdot) =  \nu (\alpha ' x + \beta 'y) \widetilde{ K} (\alpha x + \beta y ; \cdot),
\end{equation}
where
$$
\alpha + \beta = \alpha ' + \beta ' = 1,
$$
and $\nu$, $\widetilde K$ will be precised below. As explained at the end 
of Section \ref{sec:nonlocal_het}, this means that the heterogeneity in the nonlocal dispersal kernel associated with a path $(x,y)$ is no longer determined by a single point, but by the two points $\alpha x + \beta y$ and $\alpha ' x +\beta ' y$. As in the previous subsection, both points are still chosen as linear combinations of the departure and arrival points.

Notice that the diffusivity matrix associated with a dispersal kernel of the form~\eqref{ass:separation} remains the same, regardless of whether $\alpha = \alpha '$ or not, that is regardless of whether the two deciding factors are actually the same or not. Therefore one may expect that the limit local equation is also unchanged. As we will see below, it turns out that this is not the case. In some sense, the diffusivity matrix does 
not convey all the information of the motion in the heterogeneous case. This is quite similar to the observation made by two of the present authors in~\cite{KimSeo} in the kinetic model.

The motivation behind such a double heterogeneity is the fact that the dispersal kernel $K (x,y;\cdot)$ can be mostly characterized by its mass (or zero moment)
$$ m(x,y) :=  \int_{\R^N} K (x,y; z) \, dz,$$
and its relative absolute first moment
$$ g(x,y) := \frac{\displaystyle\int_{\R^N} \| z \| K (x,y;z) \, dz}{\displaystyle \int_{\R^N} K (x,y; z) \, dz}.$$
Indeed, as we explained in Section~\ref{sec:nonlocal_het}, the former can 
be interpreted as the total jump rate of individuals, and the latter as the average jump length. Both notions are crucial in the modelling of diffusion in ecology. It is then natural to allow these two functions not only to depend on space, but also to vary independently and be decided at possibly distinct points.

A typical example of dispersal kernel in the form \eqref{ass:separation} is provided by
\begin{equation}\label{typical-ex}
K (x,y;z) = \frac{m( \alpha ' x + \beta ' y)}{g( \alpha x + \beta y)^N} 
\widetilde{K} \left( \frac{z}{g(\alpha x + \beta  y)} \right),
\end{equation}
where $\widetilde{K}$ is a fixed probability density, rescaled so that \eqref{absolute} holds.

\medskip

We now proceed with the general form \eqref{ass:separation}. Dropping again the tilde for convenience, the nonlocal equation \eqref{nonlocal_het} rewrites
\begin{equation}\label{nonlocal_het_2}
\begin{array}{rcl}
u_t (t,x) & = & \displaystyle \int_{\R^N} \nu (\alpha ' y+ \beta ' x)  K (\alpha y + \beta x ; x-y)  u(t,y) \, dy \displaystyle \vspace{5pt}\\
 &  &  \quad \displaystyle  -  u (t,x) \int_{\R^N} \nu (\alpha ' x+ \beta 
' y) K (\alpha x + \beta y ; y -x)  \, dy.
\end{array}
\end{equation}
According to \eqref{ass:separation}, it is assumed that the two variables 
accounting for the heterogeneity are separable in the dispersal kernel. Moreover, in order to satisfy the assumptions of Section~\ref{sec:nonlocal_het}, we take $\nu : \R^N \to \R$ to have positive infimum and supremum, 
and $K : \R^N \times \R^N \to \R$ a nonnegative and nontrivial function. We also assume that $K(p;z)\leq \overline K(z)$ for some $\overline K \in 
L^1 (\R^N)$. As in the previous section, we impose the symmetry condition,
$$
\forall (p,z) \in \R^N  \times \R^N, \qquad K (p; -z) = K(p;z
),
$$
that both $\nu(p')$ and $K(p;z)$ are \lq\lq sufficiently smooth'', and that all derivatives of $\nu (p' ) K(p;\cdot)$ in the $p,p'$-variables are \lq\lq sufficiently $z$-integrable''.

We now consider the focusing kernels
\begin{equation}\label{Keps2}
\nu_\ep (p') := \frac{1}{\ep^2} \nu (p') \, \quad  K_\ep ( p; z):=\frac{1}{\ep^{N }}K\Big(p; \frac z \ep  \Big),\quad 0<\ep \ll 1.
\end{equation}
Plugging these into the nonlocal equation, and by a change of variable $z=\frac{y-x}{\ep}$, we obtain
\begin{eqnarray*}
\eps^2u_t(t,x)&=&\ir _{\R^N}\nu (x+\ep \alpha' z)K(x+\ep \alpha z;z)u(t,x+\ep z)\,dz\\
&&-\ir _{\R^N} \nu(x+\ep\beta'z)K(x+\ep \beta z;z)u(t,x) \,dz.
\end{eqnarray*}
After performing formal asymptotic expansions and using similar arguments 
as in subsection \ref{subsec:one-point}, we reach
\begin{eqnarray*}
u_t(t,x)&\approx & \nu(x) \ir_{\R^N}  K(x;z)\frac 12 \langle D^2u(t,x)z,z\rangle\,dz\\
&&+\alpha\nu(x) \ir_{\R^N} \Big(z\cdot \nabla_p K(x;z)\Big)\Big(z\cdot \nabla u (t,x)\Big)\,dz\\
&&+\alpha' \ir_{\R^N} (z\cdot \nabla \nu(x))K(x;z) (z\cdot \nabla u(t,x))\, dz\\
&&+\Bigg\{(2\alpha-1)\nu(x) \ir _{\R^N}\frac 12 \langle D^2_pK(x;z)z,z\rangle \,dz\\
&&\qquad + (\alpha+\alpha'-1) \ir_{\R^N} (z\cdot \nabla \nu(x))(z\cdot \nabla _p K(x;z))\,dz\\
&&\qquad+(2\alpha'-1)\ir _{\R^N}\frac 12 \langle D^2\nu (x)z,z\rangle K(x;z)\, dz\Bigg\}u(t,x).
\end{eqnarray*}

Denote
$$
q:=2-2\alpha = 2 \beta, \quad q':=2-2\alpha' = 2\beta'.
$$
In the setting \eqref{ass:separation} and after dropping the tilde, we observe that the coefficients $d_{ij}(p)$ of the diffusivity matrix $\D (p)$  of Definition~\ref{diffusivity} are given by
$$
d_{ij}(p)=\nu (p)\frac 12 \int_{\R^{N}}z_iz_jK(p;z)\,dz.
$$
In light of this, proceeding as in subsection \ref{subsec:one-point}, we can rewrite the above equation in a conservative form as
\begin{equation}\label{E2.3bis}
u_t=\sum_{i}\partial_{x_i}\Big(\sum _j\Big(\nu^{q' -q}(x)  d_{ij}^{q}(x)\partial_{x_j}\Big(\nu^{q-q'}(x)  d_{ij}^{1-q}(x)u\Big)\Big)\Big).
\end{equation}

Notice that there are two degrees of freedom, given by $q$ and~$q'$. First, find that, if $q=q'$, then \eqref{E2.3bis} is written as
\begin{equation}
\label{cas-egalite}
u_t=\sum_{i}\partial_{x_i}\Big(\sum _j\Big(d_{ij}^{q}(x)\partial_{x_j}\Big( d_{ij}^{1-q}(x)u\Big)\Big)\Big).
\end{equation}
In this case, diffusivity in the sense of Definition~\ref{diffusivity} alone decides the diffusion phenomenon. This is exactly the situation considered in  subsection \ref{subsec:one-point} with \lq\lq a single deciding 
factor''. For instance, we reach
\[
  u_t = \nabla\cdot(\nabla\cdot(\D(x) u)),\eqno(q,q')=(0,0)
\]
\[
  u_t = \nabla\cdot(\D(x)\nabla u).\ \quad\eqno(q,q')=(1,1).
\]
However, as seen in \eqref{E2.3bis}, this is no longer the case when two different types of heterogeneity are considered.

\begin{example} Let us return to our previous example \eqref{typical-ex} to better understand the meaning and the role of~$\nu$: when
$$
K (x,y;z) = \frac{m( \alpha ' x + \beta ' y)}{g( \alpha x + \beta  y)^N} \widetilde{K} \left( \frac{z}{g(\alpha x + \beta  y)} \right),
$$
then $\nu = m$ is the total jump rate. Therefore the limiting local equation \eqref{E2.3bis} typically involves not only the diffusivity matrix but also the total jump rate of the original nonlocal model. Let us also recall that, in this example, the diffusivity matrix is
$$\D (p ) = m (p) g(p)^2 \widetilde{\D},$$
where $\widetilde{\D}=(\widetilde{d_{ij}})_{ij}$ is spatially homogeneous.  Then \eqref{E2.3bis} rewrites
\begin{equation}\label{qqch}
u_t= \sum_{i}\partial_{x_i}\Big(m^{q'}(x)g^{2q}(x)\sum _j \widetilde{d_{ij}} \partial_{x_j}\Big(m^{1-q'}(x)  g^{2-2q}(x)u\Big)\Big).
\end{equation}
In this form, one may identify the role of not only the total jump rate $m$ but also the average jump length $g$.  In particular, the parameters $q$ and $(\alpha,\beta)$ correspond to the average jump length $g$ while $q'$ and $(\alpha ', \beta')$ correspond to the total jump rate $m$.
\end{example}

\begin{example}
Let us first consider the case when $\alpha '= 1$ and $\alpha =1/2$. In the above example, this means that the total jump rate is determined by the departure point, while the average jump length is determined by local conditions on the travel path, which can be reasonably simplified by choosing the middle point. From \eqref{E2.3bis}, we get
\[
  u_t = \nabla\cdot( \nu^{-1}(x) \D(x)\nabla (\nu (x) u)).\ \quad\eqno(q,q')=(1,0)
\]

Yet, we may instead consider the opposite case where the average jump length is determined at the departure point while the total jump rate is determined at the middle point (for instance, individuals may choose to go back to the departure point if they see a decline in the environmental conditions). This corresponds to $\alpha'= 1/2$ and $\alpha = 1$. From \eqref{E2.3bis}, we get
\[
  u_t = \nabla\cdot(\nu(x)\nabla\cdot(\nu^{-1}(x)\D (x) u)),\eqno(q,q')=(0,1)
\]
This equation is similar to a local equation derived in~\cite{KimSeo} from a velocity jump process. The only difference is that $\nu$ here should be $\mu^{-1}$ in~\cite{KimSeo}. To make it exactly the same, one needs to 
choose $q=0$, $q' =-1$ so that we get
\[
  u_t = \nabla\cdot(\nu^{-1}(x)\nabla\cdot(\nu(x)\D (x) u)).\eqno(q,q')=(0,-1)
\]
Still, from the point of view of the position jump process, the case $(q,q') =(0,1)$ seems more reasonable. Indeed, it is clear that diffusivity 
here is increasing with respect to $\nu$, which is typically the total jump rate. However, in the kinetic model of~\cite{KimSeo} the function $\mu$ stands for the turning frequency, with respect to which the diffusivity 
was decreasing. Thus, both notions should not be confused and it is not surprising that positive and negative exponents are swapped.
\end{example}

For general $q$ and $q'$, there are many ways to write the limiting local 
equation \eqref{E2.3bis}. We decide to write it as a \lq\lq perturbation'' of the case $(q,q')=(0,1)$ considered in the above example. To do so we write
\begin{eqnarray*}
\sum _j\nu \partial_{x_j} \Big( \nu^{-1} d_{ij}u\Big)
&=&\sum _j\nu^{q' -q}d_{ij}^q \partial_{x_j}(\nu^{q-q'}d_{ij}^{1-q}u) \\
& & \quad +\sum_j \partial_{x_j} (\nu^{q'-q -1}d_{ij}^q) \nu^{1 + q-q'}d_{ij}^{1-q}u\\
&=&\sum _j\nu^{q' -q}d_{ij}^q \partial_{x_j}(\nu^{q-q'}d_{ij}^{1-q}u) \\
& & \quad +\nu u (N_{q,q'})_i,
\end{eqnarray*}
where $N_{q,q'}(x)=((N_{q,q'})_i(x))$ is a correction vector given by
\begin{equation}\label{Eqq'}
(N_{q,q'})_i :=\sum_j  \partial_{x_j} (\nu^{q'- q -1}d_{ij}^q) \nu^{q -q'}d_{ij}^{1-q},
\end{equation}
which vanishes in the case $(q,q')=(0,1)$. As a result, the diffusion equation \eqref{E2.3bis} is recast
\begin{equation}\label{E2.7}
u_t=\nabla\cdot\Big(\nu(x)\nabla\cdot(\nu^{-1}(x)\D(x) u)-\nu(x) N_{q,q'}(x)u\Big).
\end{equation}
\begin{remark}[Orthotropic diffusion]\label{rem:ortho-2}
Under the assumption \eqref{sym-stonger} for orthotropic diffusion, equation \eqref{E2.3bis} turns into
\[
u_t=\sum_{i}\partial_{x_i}\Big(\nu^{q' -q} (x)d_{ii}^{q}(x)\partial_{x_i} \Big(\nu^{q-q'}(x) d_{ii}^{1-q}(x)u\Big)\Big).
\]
Furthermore, since the diffusivity matrix $\D(p)$ is diagonal, $\D^q(p)$ can be defined element-wise. Then, we may write it as
\begin{equation}\label{Eortho}
u_t=\nabla\cdot\Big(\nu^{q'-q}(x)\D^{q}(x) \Big(\nabla \cdot (\nu^{q-q'}(x)\D^{1-q}(x)u)\Big)\Big).
\end{equation}
\end{remark}

\begin{remark}[Isotropic diffusion]\label{rem:iso-2}
If $z\mapsto K(p;z)$ is radial symmetric, we obtain
$$
\D(p)=D(p)\mathbf{I},\quad D(p):=\frac{\nu(p)}{2N} \ir _{\R^N} \|  z\| ^2 K(p;z)\,dz.
$$
The isotropic diffusion equation is now
\begin{equation}\label{Eiso}
u_t=\nabla\cdot\Big(\nu^{q'-q}(x)D^{q}(x) \nabla \Big(\nu^{q-q'}(x)D^{1-q}(x)u\Big)\Big).
\end{equation}
For instance, if $(q,q')=({1\over2},1)$, we obtain
\[
u_t=\nabla\cdot\Big(\sqrt{\nu(x) D(x)}\, \nabla \big(\sqrt{\nu^{-1}(x)D(x)}\, u\big)\Big).
\]
\end{remark}

\section{Steady states and numerical simulations}\label{sec:steady}

In this section we present some  enlightening
observations  regarding the steady states of \eqref{nonlocal_het_gen},
when the jump rate is of the form
\begin{equation}\label{ex1-nouveau}
J(x,y)=  m(\alpha x+\beta y) \widetilde{K} (y-x).
\end{equation}
This connects with the issue discussed  in, e.g., \cite{Cos-Dav-Mar-12, Can-Cos-Lou-Rya-12} and is complemented by some numerical simulations. 
Those simulations are performed
on a bounded domain,  $B(0,R)\subset\mathbb{R}^{N}$,
see subsection \ref{subsec:NumImp} for details.

\subsection{Steady states}\label{subsec:Steady-states}

We consider the case of a single deciding factor with $J(x,y)$ of
the form \eqref{ex1-nouveau}. In what follows, we omit the tilde 
notation for readability. The nonlocal model \eqref{nonlocal_het_gen}
is then recast
\begin{equation}
u_{t}=\int_{\mathbb{R}^{N}}\left[m(\alpha y+\beta x)u(t,y)-m(\alpha x+\beta y)u(t,x)\right]K(y-x)dy,\quad\alpha+\beta=1.\label{eq:nonlocal}
\end{equation}
As seen in subsection \ref{subsec:one-point}, taking the focusing
kernel limit leads, in the isotropic case, to the local model 
\begin{equation}
u_{t}=\nabla\cdot\left(D^{q}(x)\nabla\left(D^{1-q}(x)u\right)\right),\qquad q=2-2\alpha,\label{eq:local}
\end{equation}
with $D(x)=\frac{1}{2}m(x) k$,  where $k:=\frac 1N \int_{\R^N}\Vert z\Vert ^2 K(z)dz$. For the rest of this section, we consider
a nonnegative, nontrivial initial data $u_{0}\in L^{1}(\mathbb{R}^{N})$.
Also, we assume that $m$ is a nonnegative continuous function. Finally, for numerical simulations and unless otherwise stated, we shall use the Gaussian kernel
$$
K(z):=C_N e^{-a_N\Vert z\Vert ^{2}},
$$
where $C_N>0$, $a_N>0$ are such that not only $\int_{\R^{N}}K(z)dz=1$ but also $\int_{\R^{N}}\Vert z\Vert K(z)dz=1$. In particular $C_1=a_1=\frac 1 \pi$.\vspace{5pt}

\noindent \textbf{The case $\alpha=0$ ($q=2$).} In this case, we observe  that $u(x)\equiv Cm(x)$ is a steady state of \eqref{eq:nonlocal}
for any constant $C\geq0$. This is consistent with the local counterpart,
that is \eqref{eq:local}, for which $u(x)\equiv CD(x)$ is also a steady
state. In other words, at equilibrium, the population is mainly localized
at positions~$x$ where $m(x)$, or equivalently $D(x)$, is large.
This is in agreement with the model since the total jump rate from $x$ to 
$y$ is given by
$m(y)$, meaning the jump rate is higher for jumps
that arrive at high values of $m$, and conversely smaller for jumps
that arrive at low values of $m$. Therefore the solution of \eqref{eq:nonlocal}
is expected to converge, at large time, to
\begin{equation}
p(x):= Cm(x),\qquad C=\frac{\int_{\mathbb{R}^{N}}u_{0}(x)dx}{\int_{\mathbb{R}^{N}}m(x)dx}\geq0,\label{eq:p_alpha0}
\end{equation}
the value of $C$ ensuring the mass conservation. In particular, if
$m\notin L^{1}(\mathbb{R}^{N})$, we expect that $p(x)\equiv0$.

\begin{figure}
\centerline{\includegraphics[scale=0.45]{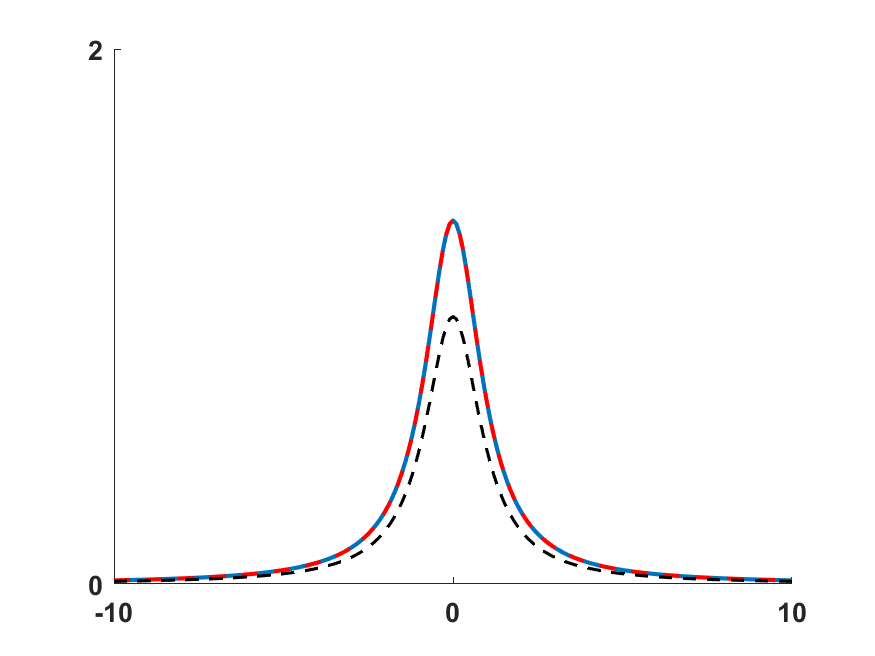}\includegraphics[scale=0.45]{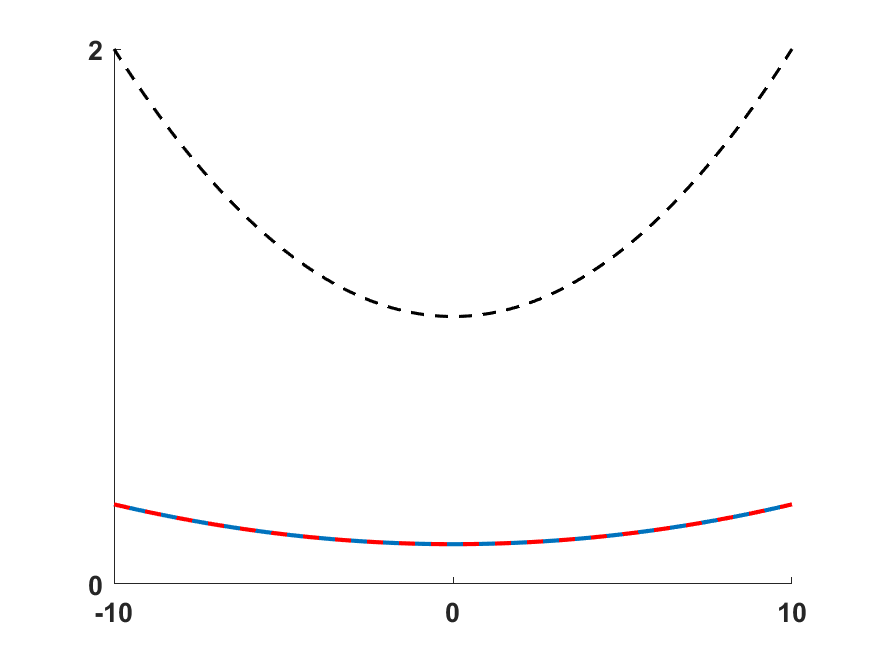}}

\caption{Here $\alpha=0$ and $R=10$. In dotted black, the function $m$.
In blue, the numerical approximation of $u$ at time $t=4000$. In dotted
red, the expected profile $p_{R}(x)$. On the left $m(x)=\frac{1}{1+x^{2}}$,
and on the right $m(x)=1+\frac{1}{100}x^{2}$. The outcome depends only on the initial mass, here $\int_{\Omega_{R}}u_{0}(x)dx=4$.\label{fig:alpha0}}
\end{figure}

We now present some numerical simulations to illustrate the above.
Since these simulations are done on the bounded domain $\Omega_{R}=B(0,R)\subset\mathbb{R}^{N}$,
the expected profile is (see subsection \ref{subsec:NumImp} for details)
\begin{equation}
p_{R}(x):= C_R m(x),\qquad C_R=\frac{\int_{\Omega_{R}}u_{0}(x)dx}{\int_{\Omega_{R}}m(x)dx}>0,\label{eq:p_alpha0_R}
\end{equation}
which is consistent with our simulations, see Figure \ref{fig:alpha0}. Moreover, we readily check that $p_{R}\to p$ uniformly on compact
sets of $\mathbb{R}^{N}$, and uniformly on $\mathbb{R}^{N}$ for
bounded $m$. However, for general $m$, this convergence may not
be uniform on $\mathbb{R}^{N}$, as can be seen with $N=1$, $m(x)=e^{x}$,
for which $||p_{R}||_{\infty}\to \int_{\mathbb{R}}u_{0}(x)dx$ while $p\equiv0$. \vspace{5pt}

\noindent \textbf{The case $\alpha=1$ ($q=0$).} Here, we assume $m(x)>0$ for all $x$. We observe that $u(x)\equiv Cm(x)^{-1}$
is a steady state of \eqref{eq:nonlocal} for any constant $C\geq0$.
This is again consistent with the local counterpart, that is \eqref{eq:local},
for which $u\equiv CD(x)^{-1}$ is also a steady state. The situation
is thus reversed compared to the case $\alpha=0$: the population
will mainly concentrate at positions $x$ where $m(x)$ is small. This is,
again, in agreement with the model since the jump rate from $x$ to $y$ is 
given by
$m(x)$, meaning the jump rate is higher for jumps
that depart from high values of $m$, and conversely smaller for jumps
that depart from low values of $m$. Therefore the solution
of \eqref{eq:nonlocal} is expected to converge, at large time, to
\begin{equation}
p(x):= Cm(x)^{-1},\qquad C=\frac{\int_{\mathbb{R}^{N}}u_{0}(x)dx}{\int_{\mathbb{R}^{N}}\frac{1}{m(x)}dx}\geq0,\label{eq:p_alpha1}
\end{equation}
the value of $C$ ensuring the mass conservation. In particular, if
$m^{-1}\notin L^{1}(\mathbb{R}^{N})$, we expect that $p(x)\equiv0$. 

\begin{figure}
\centerline{\includegraphics[scale=0.45]{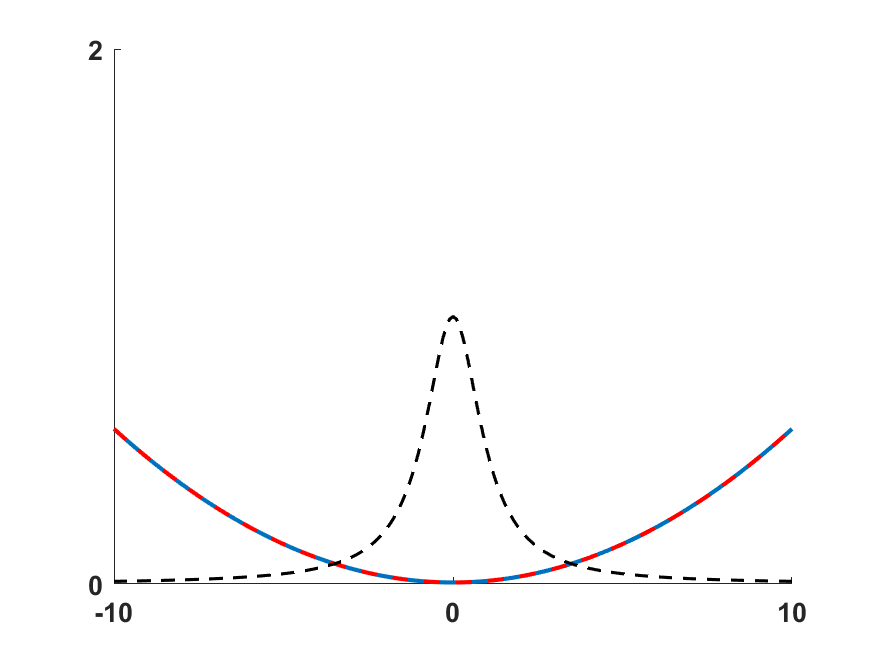}\includegraphics[scale=0.45]{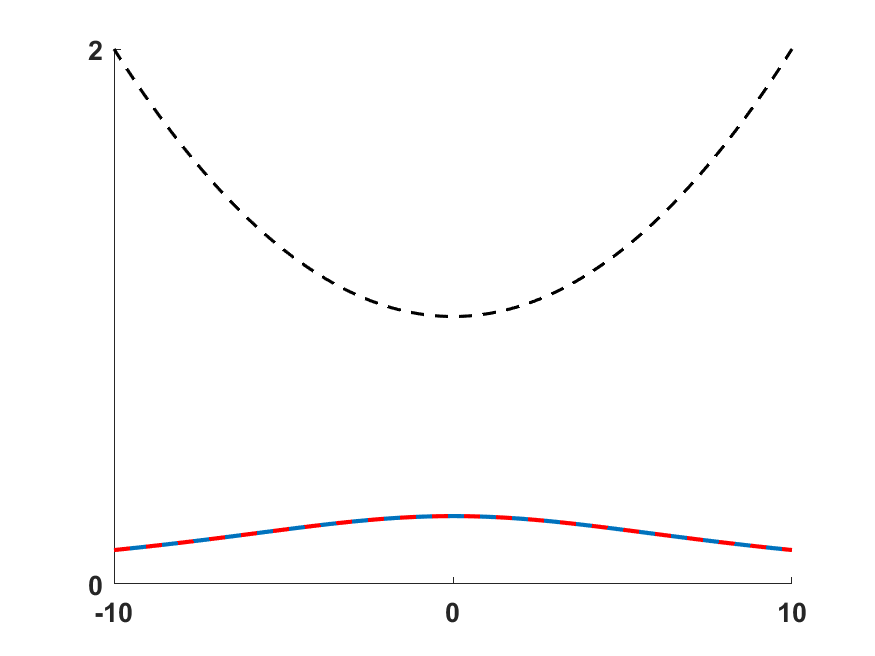}}

\caption{Conditions are the same as in Figure \ref{fig:alpha0}, except $\alpha=1$.\label{fig:alpha1}}
\end{figure}

Similarly to the case $\alpha=0$, on the bounded domain $\Omega_{R}=B(0,R)\subset\mathbb{R}^{N}$,
the expected profile would be
\begin{equation}
p_{R}(x):= C_R m(x)^{-1},\qquad C_R=\frac{\int_{\Omega_{R}}u_{0}(x)dx}{\int_{\Omega_{R}}\frac{1}{m(x)}dx}>0,\label{eq:p_alpha1_R}
\end{equation}
which is consistent with our simulations, see Figure \ref{fig:alpha1}.
As in the case $\alpha=0$, we readily check that $p_{R}\to p$ uniformly 
on compact
sets of $\mathbb{R}^{N}$.\vspace{5pt}

\begin{figure}
\centerline{\includegraphics[scale=0.45]{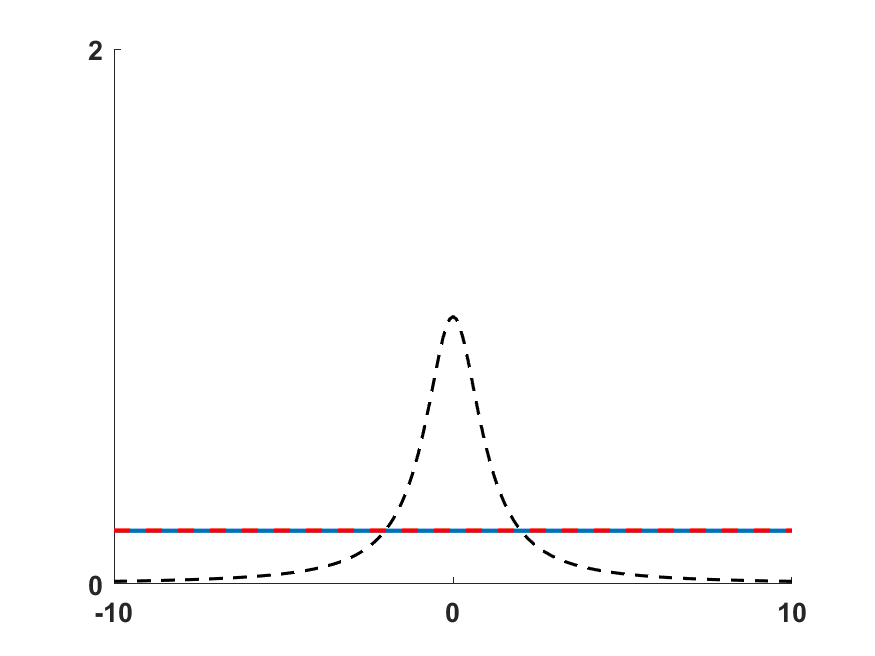}\includegraphics[scale=0.45]{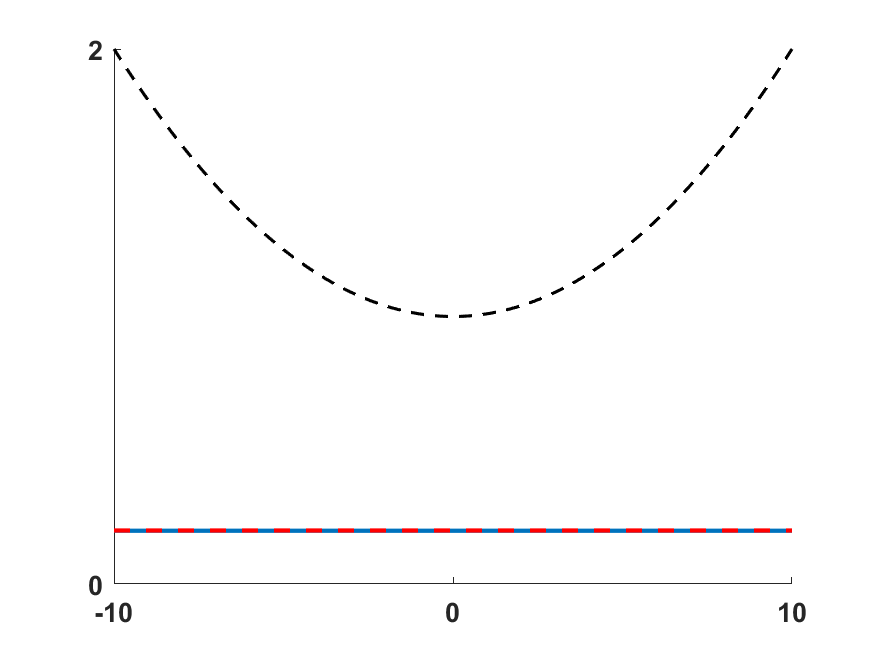}}

\caption{Conditions are the same as in Figure \ref{fig:alpha0}, except $\alpha=\frac{1}{2}$.\label{fig:alpha05}}
\end{figure}

\noindent \textbf{The case $\alpha=\frac{1}{2}$ ($q=1$).} In that setting, it is clear that constants are steady states of both
\eqref{eq:nonlocal} and \eqref{eq:local}. Therefore the profile
of $m$ plays no role in the distribution of population at equilibrium.
This is again expected since then
$$
J(x,y)=m\left(\frac{x+y}{2}\right)K(y-x)=J(y,x),
$$
that is the jump rate from $x$ to $y$ is equal to that from $y$ to $x$. Therefore the solution of \eqref{eq:nonlocal}
is expected to converge, at large time, to $p(x)\equiv0$, due to the mass 
conservation.
Note that on a bounded domain $\Omega_{R}=B(0,R)\subset\mathbb{R}^{N}$,
we however expect a convergence towards the constant $p_{R}=\frac 1{\left|\Omega_{R}\right|}\int_{\Omega_{R}}u_{0}(x)dx$,
which is in agreement with our numerics, see Figure \ref{fig:alpha05}. Clearly, $p_R\to 0$ as $R\to +\infty$.\vspace{5pt}

\noindent \textbf{Towards general $0\protect\leq\alpha\protect\leq1$.} The steady states are thus well-understood for general functions $m$
in the cases $\alpha\in\{0,\frac{1}{2},1\}$. In particular, those
steady states depend only on the total jump rate $m$, and not on
the jump distribution $K$. We aim at extending those results for
general $\alpha$. Hence, for a given $\alpha\not\in\{0,\frac{1}{2},1\}$ and a given nonnegative and nontrivial  $m$, we look for a nonnegative and 
nontrivial $u$ such
\begin{equation}
m(\alpha y+\beta x)u(y)=m(\alpha x+\beta y)u(x),\qquad\forall x,y\in\mathbb{R}^{N},\label{eq:um_rela}
\end{equation}
which makes $u=u(x)$ a steady state for \eqref{eq:nonlocal}. Notice that such a function~$u$ is not ensured to exist, as can be seen with
$m(\cdot)=||\cdot||$: in this case, for $x=0$, \eqref{eq:um_rela} enforces $u(y)\equiv\frac{\beta}{\alpha}u(0)$
for any $y\neq0$; returning to \eqref{eq:um_rela}, this is absurd  if $\alpha\neq\frac{1}{2}$. 

\begin{figure}
\centerline{\includegraphics[scale=0.5]{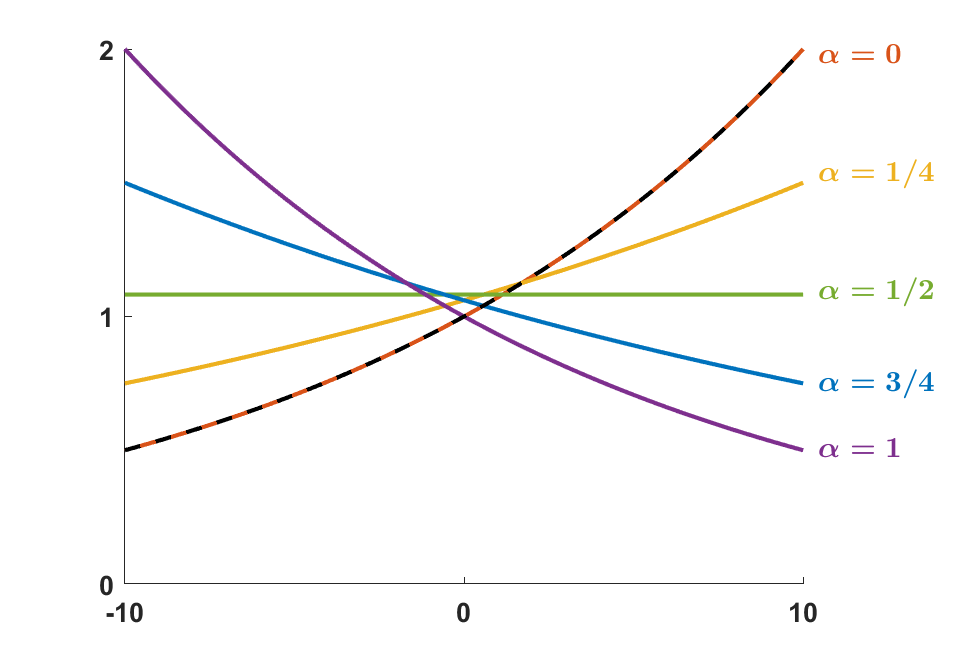}}
\caption{In dotted black, the function $u_{0}(x)=m(x)=e^{ax}$ with $a=\frac{\ln2}{10}$. In solid lines, for various $\alpha$, the numerical approximation of $u$ at $t=4000$, which all coincide with $p_{R}^{\alpha}$. \label{fig:expoTrans}}
\end{figure}

However, there are some profiles $m$ for which such a function $u$
exists for any $0\leq\alpha\leq1$. Indeed, assuming $N=1$, $m(x)=e^{ax}$
with $a\in\mathbb{R}$, we see that for each $\alpha$, the function
$u^{\alpha}(x)=e^{(\beta-\alpha)ax}=m(x)^{\beta-\alpha}$ solves \eqref{eq:um_rela}.
Because $u^{\alpha}\notin L^{1}(\mathbb{R})$, we expect the solution
of \eqref{eq:nonlocal} to converge towards $p^{\alpha}(x)\equiv0$
due to the mass conservation. However, on the bounded domain $\Omega_{R}=[-R,R]$,
the expected profile is
\begin{equation}
p_{R}^{\alpha}(x):=C^{\alpha}_R u^{\alpha}(x),\qquad C^{\alpha}_R=\frac{\int_{\Omega_{R}}u_{0}(x)dx}{\int_{\Omega_{R}}u^{\alpha}(x)dx}>0,\label{eq:p_alpha_R}
\end{equation}
which is consistent with our simulations, see Figure \ref{fig:expoTrans}, for which we have chosen $u_0\equiv m$ which, obviously, is not in $L^1(\R)$. Nevertheless, our goal here is to illustrate the fact that, as $\alpha$ increases, there is a smooth transition from~$m$ (when $\alpha=0$) to $m^{-1}$ (when $\alpha=1$), the \lq\lq switch'' precisely occurring at $\alpha=\frac 12$.  Last, one can extend this example to $N\geq 2$ with $m(x)=e^{\langle a,x\rangle}$
and $a\in\mathbb{R}^{N}$, for which $u(x)=e^{(\beta-\alpha)\langle a,x\rangle}$
solves \eqref{eq:um_rela}.

\begin{figure}
\centerline{\includegraphics[scale=0.5]{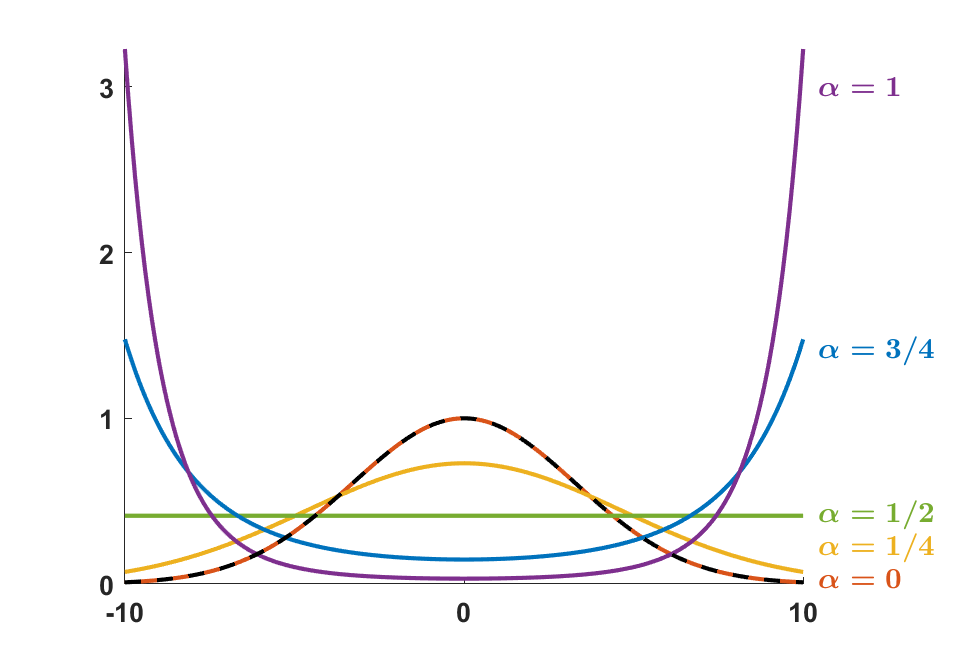}}
\caption{In dotted black, the function $u_{0}(x)=m(x)=e^{-ax^2}$ with 
$a=\frac{\ln100}{100}$. In solid lines, for various $\alpha$, the numerical approximation of $u$ at $t=4000$, which all coincide with $p_{R}^{\alpha}$.\label{fig:GaussTrans}}
\end{figure}

Another example with $N=1$ is given by the Gaussian profile $m(x)=e^{-ax^{2}}$
with $a>0$, for which $u^{\alpha}(x)=e^{-(\beta^{2}-\alpha^{2})ax^{2}}$
solves \eqref{eq:um_rela}. Since $u^{\alpha}\in L^{1}(\mathbb{R})$
if and only if $\alpha<\frac{1}{2}$, we expect that the solution
of \eqref{eq:nonlocal} converges towards
\[
p^{\alpha}(x)=\begin{cases}
C^{\alpha}u^{\alpha}(x) & \text{if }\alpha<\frac{1}{2},\\
0 & \text{if }\alpha\geq\frac{1}{2},
\end{cases}\qquad C^{\alpha}=\frac{\int_{\mathbb{R}^{N}}u_{0}(x)dx}{\int_{\mathbb{R}^{N}}u^{\alpha}(x)dx}>0.
\]
Meanwhile, on the bounded domain $\Omega_{R}=[-R,R]$, since $u^{\alpha}\in L^{1}(\Omega_{R})$
for all $\alpha$, we expect that the solution converges towards $p_{R}^{\alpha}$
given by \eqref{eq:p_alpha_R}. This is what we observe numerically,
see Figure \ref{fig:GaussTrans}. Note that the above comments for
the exponential profile still apply here. Last, one can extend
this example to $N\geq 2$ with $m(x)=e^{-\langle a,x^{2}\rangle}$,
where $a\in(0,+\infty)^{N}$ and $x^{2}:=(x_{1}^{2},\dots,x_{N}^{2})$,
for which $u(x)=e^{-(\beta^{2}-\alpha^{2})\langle a,x^{2}\rangle}$
solves \eqref{eq:um_rela}.\vspace{5pt}

\begin{figure}
\centerline{\includegraphics[scale=0.5]{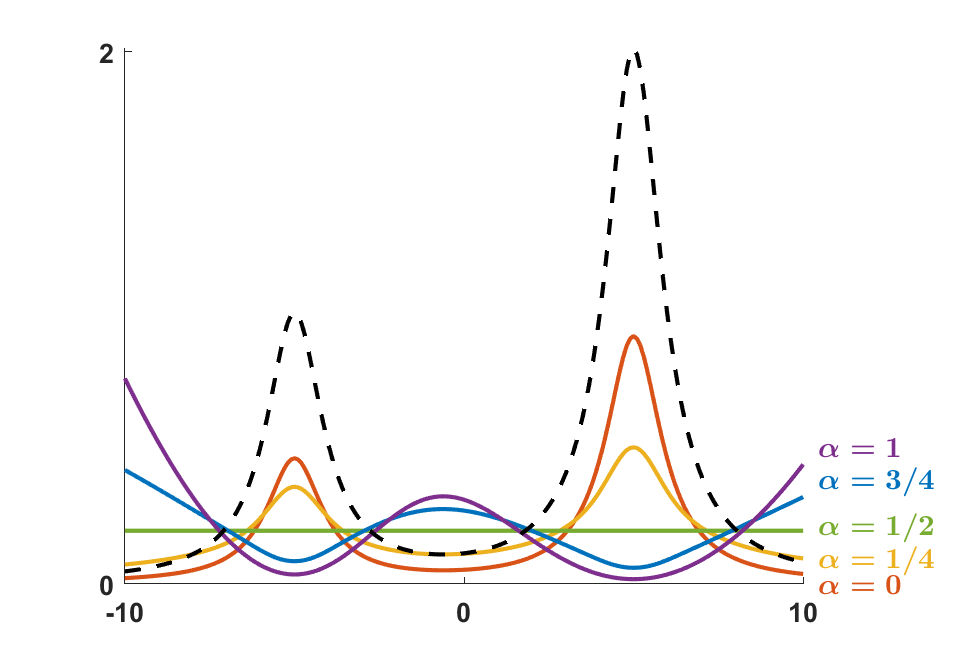}}\caption{In dotted black, the function $m$ given by \eqref{eq:m_asym}. In solid lines, for various $\alpha$, the numerical approximation of $u$ at $t=4000$ . The initial data is of mass $\int_{\Omega_{R}}u_{0}(x)dx=4$.\label{fig:asym_m}}
\end{figure}

\noindent \textbf{Where are the individuals? The role of $\alpha$.} Our above observations and results can be summarized as follows: when
$\alpha=0$ (resp. $\alpha=1$, resp. $\alpha=\frac{1}{2}$), some
steady states are proportional to $m$ (resp. $m^{-1}$, resp. $x\mapsto1$).
We therefore conjecture that, as $\alpha$ increases, the population
should concentrate at positions $x$ where $m(x)$ is small. This 
was proved above for the particular cases of exponential and Gaussian profiles of $m$.
To further support this conjecture, we investigate the case of a two-patch asymmetric function $m$, namely
\begin{equation}
m(x)=\gamma(x+5)+2\gamma(x-5),\qquad\gamma(x)=\frac{1}{1+x^{2}}.\label{eq:m_asym}
\end{equation}
The numerical solutions are displayed in Figure \ref{fig:asym_m},
for different values of $\alpha$. We observe indeed that, as $\alpha$
increases, the population shifts away from the positions where  $m$ is large and settles in positions where  $m$ is small.\vspace{5pt}

\noindent \textbf{The role of the tails of $K$.} In all the above simulations, we always considered that the probability
density $K$ (rescaled so that $\int_{\mathbb{\mathbb{R}}^{N}}||z||K(z)dz=1$)
has Gaussian tail, that is (say $N=1$)
$$
K(z)=\frac{1}{\pi}e^{-\frac{1}{\pi}z^{2}}=:K_1(z).
$$
In order to draw a comparison, we also consider the following kernels:
\[
K_{2}(z):=\frac{1}{2}e^{-|z|},\qquad K_{3}(z):=\frac{1}{\pi}\frac{1}{1+z^{4}/4}.
\]
Note that, denoting $k_{i}=\int_{\mathbb{R}}z^{2}K(z)dz$,
we have $k_{1}=\frac{\pi}{2}$ and $k_{2}=k_{3}=2$.
Finally, we consider one last kernel $K_4$ with a heavier tail, namely
\[
K_4(z):= \frac{1}{1+|z|^3} \left( \frac{A}{\log(2+|z|)^3} 
- \frac{B}{\log(2+|z|)^5} +\frac{C}{\log(2+|z|)^7} \right),
\]
where $A\approx 2.32$, $B\approx 2.38$ and  $C\approx 0.62$ are numerically selected such that $\int_\R K_4(z)dz=1$,  $\int_{\mathbb{\mathbb{R}}}|z|K_4(z)dz=1$ and $k_4=2$. Notice that $K_4\geq 0$ does hold and that $K_4$ is not radially decreasing. In particular, for $i\geq 2$, all the kernels $K_i$ have the same second moment, which means that they all lead to the same local model in the focusing kernel limit. 

\begin{figure}
\includegraphics[scale=0.45]{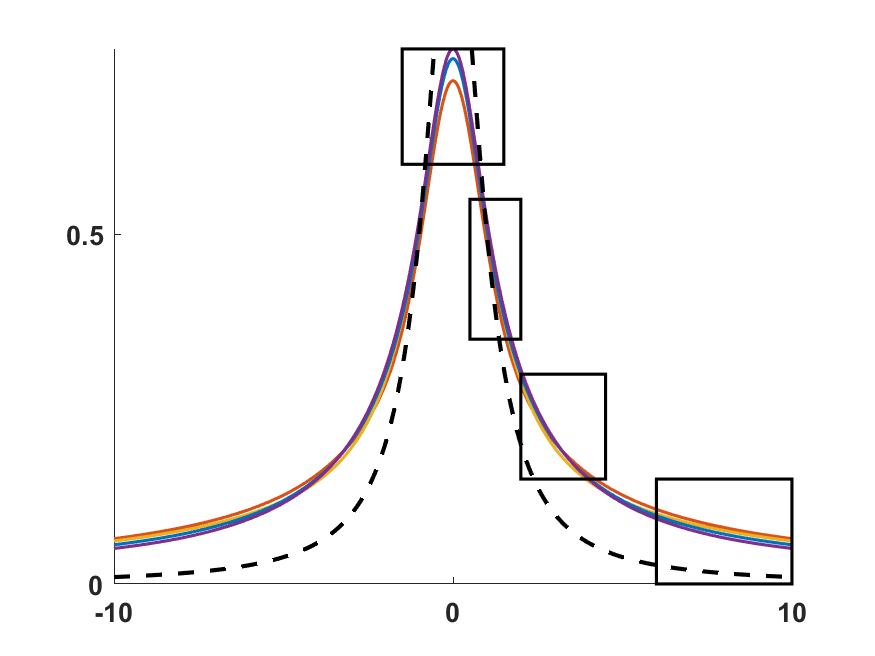}\\
\centerline{\includegraphics[scale=0.45]{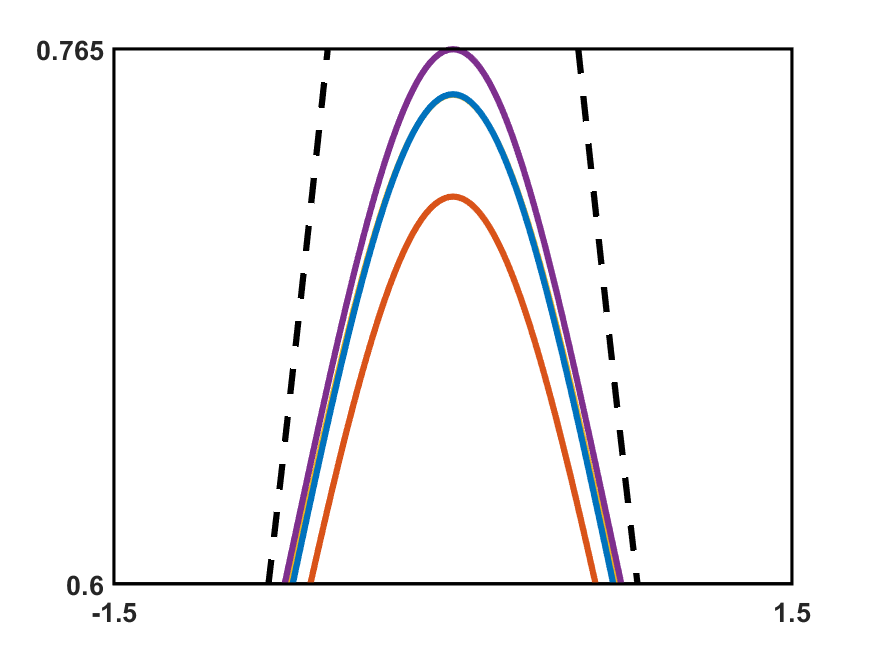}
\includegraphics[scale=0.45]{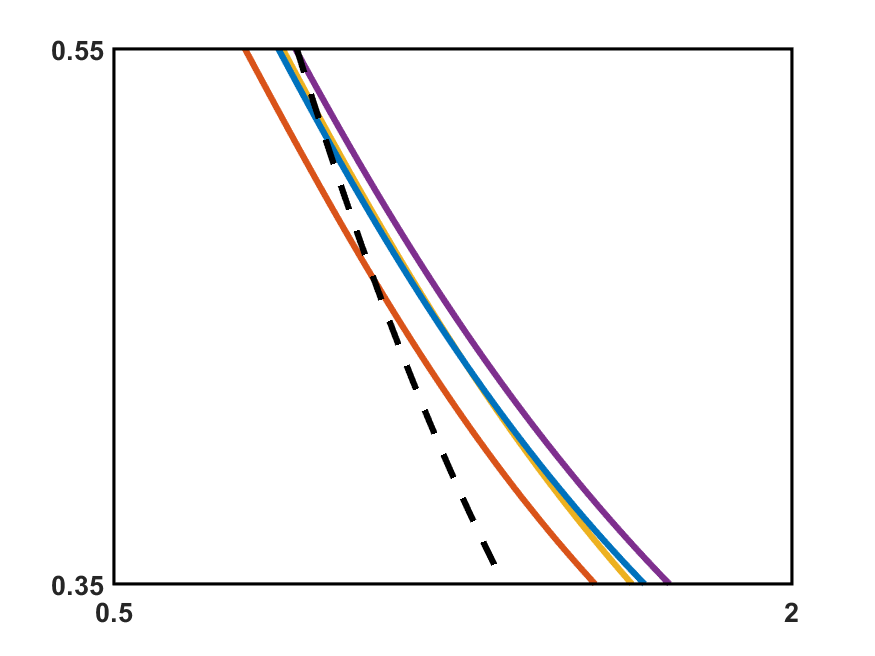}}
\centerline{
\includegraphics[scale=0.45]{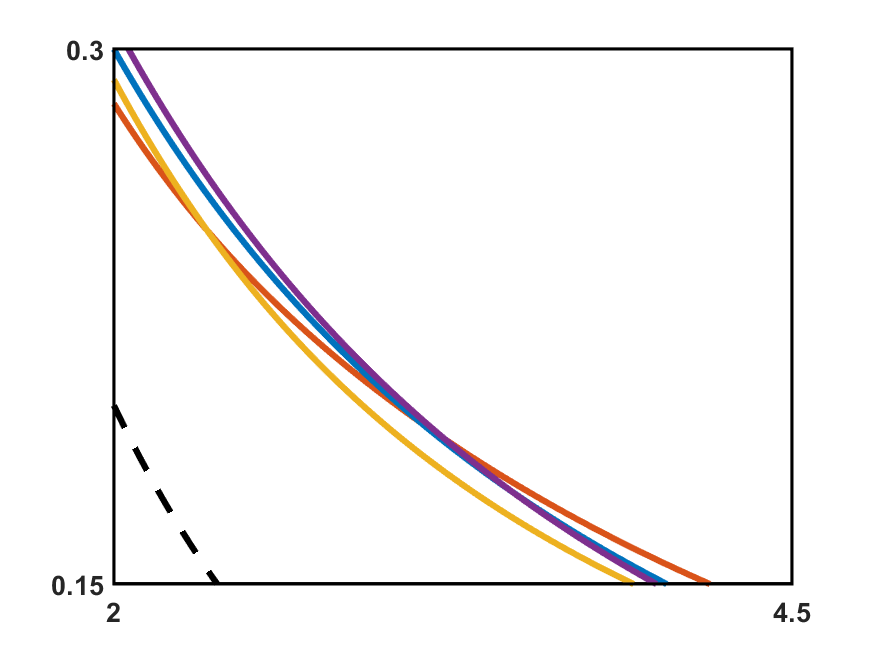}
\includegraphics[scale=0.45]{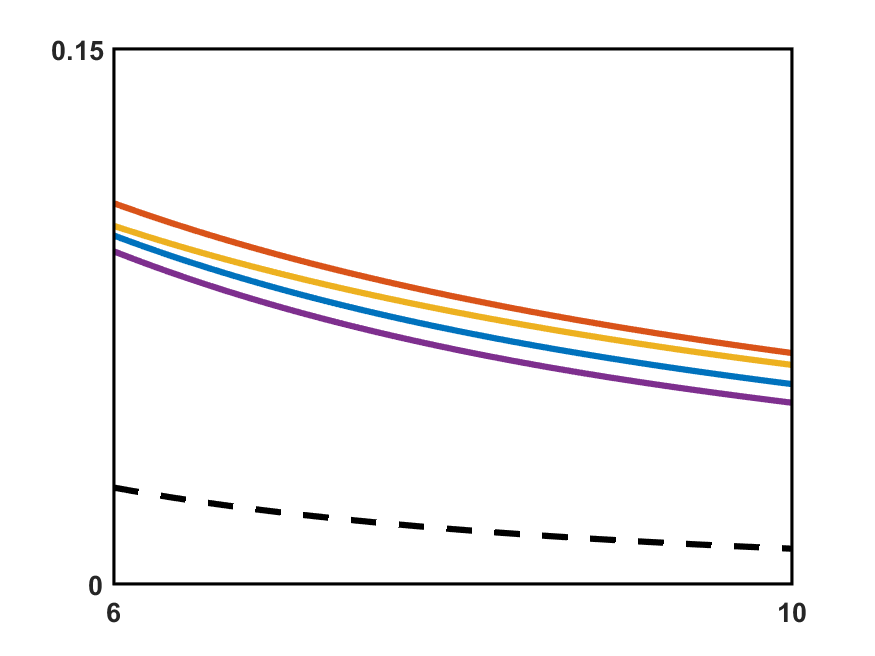}}
\caption{Here $\alpha=\frac{1}{4}$. In dotted black, the function $m(x)=\frac{1}{1+x^{2}}$. In red, the
numerical approximation of the solution at $t=10^4$, corresponding to the kernel $K_{1}$ (red),  $K_{2}$ (yellow), $K_{3}$ (blue) and  $K_4$ (purple). Some parts of the first image have been zoomed in to improve visibility. On the first and second zooms, the yellow curve almost coincides with the blue one. The initial data is of mass $\int_{\Omega_{R}}u_{0}(x)dx=4$.
The results are similar on larger domains.\label{fig:Kvariable}}
\end{figure}

\begin{figure}
\centerline{
\includegraphics[scale=0.45]{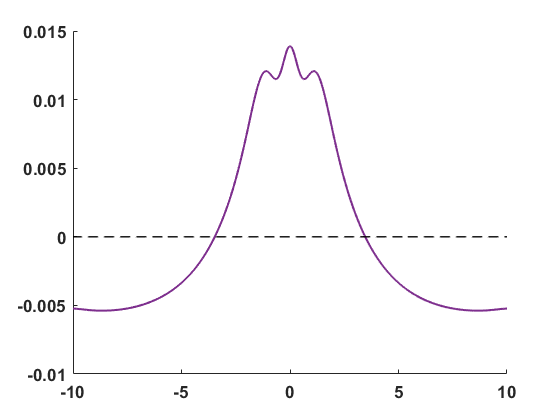}
\includegraphics[scale=0.45]{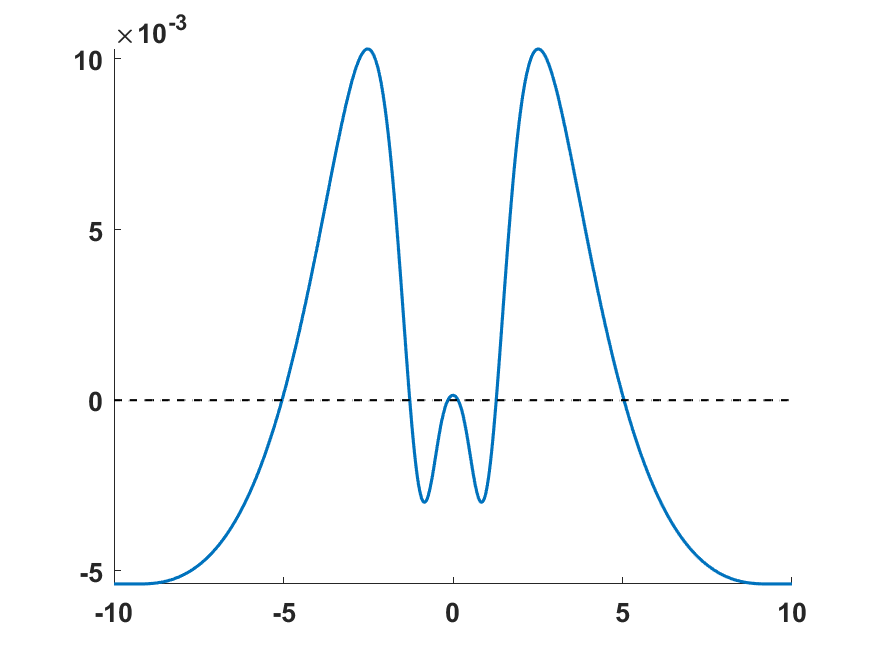}
}
\caption{Left: $u_4-u_3$. Right: $u_3-u_2$.  Here $u_i$ corresponds to the numerical solution when $K=K_i$, in the conditions of Figure \ref{fig:Kvariable}.\label{fig:uj-ui}}
\end{figure}

Our goal here is to numerically compare the steady states  depending on the chosen kernel. When $\alpha\in\{0,\frac{1}{2},1\}$, we do not
observe, as expected, any difference, for various functions $m$. The situation is quite different when, say, $\alpha=\frac{1}{4}$, see Figure \ref{fig:Kvariable}. Indeed, we observe the following: the heavier the tails of $K_i$, the more the corresponding solution  $u_i$ tends to concentrate at positions where
$m$ is large, as can be seen from the first and last zoomed parts. 
Knowing that, one may expect that $u_j-u_i$ changes sign only once on $(0,+\infty)$. It turns out to be the case for $(i,j)=(4,3)$ but not for $(i,j)=(2,3)$, see Figure \ref{fig:uj-ui}. The latter also indicates that, in a neighbourhood of $x=0$,  $u_3- u_2$ is (positive and) very small while  $u_4$ and $u_3$ are clearly distinct.

\subsection{Numerical implementation\label{subsec:NumImp}}

Since any numerical computation has to be done on a bounded domain
$\Omega\subset\mathbb{R}^{N}$, the numerical integration of the right-hand
side of \eqref{nonlocal_het_gen} is not trivial, even for
simple cases such as $J(x,y)=K(y-x)$. Indeed, in order to compute
the convolution $K*u$ close to the boundary $\partial\Omega$,
one would have to extrapolate $u$ outside of $\Omega$ itself. To
circumvent this, we consider the following problem on the bounded domain $\Omega\subset \R^N$:
\begin{equation}
u_{t}=\int_{\Omega}\left[m(\alpha y+\beta x)u(t,y)-m(\alpha x+\beta y)u(t,x)\right]K(y-x)dy.\label{eq:nonlocal_bounded}
\end{equation}
In other words, we assume that any jump that would leave $\Omega$
is in fact omitted. The functions $p_{R}$ and $p_{R}^{\alpha}$
appearing in subsection \ref{subsec:Steady-states} are in fact steady
states of \eqref{eq:nonlocal_bounded} with $\Omega=\Omega_{R}=B(0,R)\subset\mathbb{R}^{N}$, and we have checked that our numerical solution converge towards them for some large values of~$R$.
As mentioned in subsection \ref{subsec:Steady-states}, the steady
states $p_{R}$ and $p_{R}^{\alpha}$ typically converge, as $R\to+\infty$, 
to
a steady state of \eqref{eq:nonlocal}, locally uniformly in $\mathbb{R}^{N}$.

\section{Summary and perspectives}\label{s:summarize}

For the modelling of multi-dimensional nonlocal heterogeneous diffusion, we have proceeded as follows: to any possible path $(x,y)$ we have assigned a dispersal kernel $K(x,y;\cdot)$ such that $K(x,y;-z)=K(x,y;z)$, and relevant quantities for this path are
$$
m(x,y):=\int_{\R^N} K(x,y;z)\,dz \quad \text{ the total jump rate,}
$$
and
$$
g(x,y):=\frac{\displaystyle{\int_{\R^N} \|z \|  K(x,y;z)\, dz}}{\displaystyle{\int _{\R^{N}} K(x,y;z)\,dz}} \quad \text{ the average jump length.}
$$
Also, for $p\in \R ^N$, we have defined
$$
\D(p):=\left(d_{ij}(p):=\frac {1}{2}\ir _{\R^N}   z_iz_j K(p,p;z) \, dz\right)_{1\leq i,j\leq n} \quad \text{ the diffusivity matrix.}
$$

\subsection{A single deciding factor}\label{ss:one-DF} First we have considered the case
$$
K(x,y;z)\leftarrow K(\alpha x+\beta y;z), \quad \alpha+\beta=1.
$$
We have revealed, in the focusing kernel limit, the connection
$$
\fbox{\text{ Nonlocal model  \eqref{nonlocal_het_1} $\xrightarrow{q=2-2\alpha}$  diffusion equation  \eqref{componentwise},  or  \eqref{ANq}},}
$$
which is possibly anisotropic. In the so-called orthotropic case, the diffusivity matrix is diagonal and the limit diffusion equation is recast
\[
u_t=\nabla\cdot\Big(\D^{q}(x) \Big(\nabla \cdot (\D^{1-q}(x)u)\Big)\Big).
\]
In the so-called isotropic case, the diffusivity matrix is scalar and the 
limit diffusion equation is recast
\[
u_t=\nabla\cdot\Big(D^{q}(x) \nabla \Big(D^{1-q}(x)u\Big)\Big).
\]
The following examples recover some standard diffusion laws and are worth 
being mentioned.

\begin{center}
\begingroup
\setlength{\tabcolsep}{5pt} % Default value: 6pt
\renewcommand{\arraystretch}{1.5} % Default value: 1
\begin{tabular}{|l|l|l|}
  \hline
Deciding factor & Limit diffusion equation  & Limit  law \\
  \hline
$\alpha=1$, departure point & $u_t=\Delta (D(x)u)$  & Chapman's law\\
$\alpha=\frac 34$, \lq\lq close to departure''  & $u_t=\nabla\cdot\left(\sqrt{D(x)} \nabla (\sqrt{D(x)} u)\right)$  & Wereide's law\\
$\alpha=\frac 12$, middle point  & $u_t=\nabla\cdot\left(D(x) \nabla u\right)$  & Fick's law\\
$\alpha=0$, arrival point  & $u_t=\nabla\cdot\left(D^{2}(x) \nabla \left(\frac u{D(x)}\right)\right)$  & \\
  \hline
\end{tabular}
\endgroup
\end{center}\vspace{5pt}

Furthermore, in the framework $K(\alpha x+\beta y;z)\leftarrow m(\alpha x+\beta y)K(z)$, we have investigated the form of the steady state solutions. This revealed that, as $\alpha$ increases, the population shifts away 
from the positions where $m$ is large and settles in positions where $m$ is small. In some situations, the steady states are fully determined by $\alpha$ and $m$. Nevertheless, when $\alpha \notin \{0,\frac 12,1\}$, we have tested the influence of different kernels $K$ and noticed some subtle differences on the shape of the steady states. It seems that, the heavier the tails of~$K$, the more the solution concentrates, which is slightly counter-intuitive. Understanding how the steady states are determined by the interplay of $\alpha$, $m$ and $K$ appears very challenging and deserves further investigations.

\subsection{Two deciding factors}\label{ss:two-DF} Next we have considered the case
$$
K(x,y;z)\leftarrow \nu (\alpha'x+\beta ' y)K(\alpha x+\beta y;z), \quad \alpha+\beta= \alpha'+\beta'=1.
$$
We have revealed, in the focusing kernel limit, the connection
$$
\fbox{\text{ Nonlocal model  \eqref{nonlocal_het_2} $\xrightarrow[q'=2-2\alpha ']{q=2-2\alpha}$  diffusion equation  \eqref{E2.3bis},  or  \eqref{E2.7}},}
$$
which is possibly anisotropic. In the so-called orthotropic case, we refer to Remark \ref{rem:ortho-2}.
In the so-called isotropic case, we refer to Remark \ref{rem:iso-2} and the following examples are worth being mentioned (where D.F. means deciding factor).

\begin{center}
\begingroup
\setlength{\tabcolsep}{5pt} % Default value: 6pt
\renewcommand{\arraystretch}{1.5} % Default value: 1
\begin{tabular}{|l|l|l|}
  \hline
 D.F. for average jump length & D.F. for total jump rate & Limit diffusion equation  \\
  \hline
$\alpha=1$, departure point & $\alpha'=\frac 12$, middle point &$ u_t 
= \nabla\cdot(\nu(x)\nabla\cdot(\nu^{-1}(x)D (x) u))$ \\
$\alpha=1$, departure point  & $\alpha'=\frac 32$, \lq\lq strange'' point & $u_t=\nabla\cdot(\nu^{-1}(x)\nabla\cdot(\nu(x)D (x) u))$  \\
$\alpha=\frac 12$, middle point  & $\alpha'=1$, departure point &  $u_t=\nabla\cdot(\nu^{-1}(x)D(x)\nabla\cdot(\nu(x) u))$\\
$\alpha=\frac 34$, \lq\lq close to departure''  & $\alpha'=\frac 12$, 
middle point &  $u_t=\nabla\cdot(\sqrt{\nu D(x)}\nabla\cdot(\sqrt{\nu^{-1} D(x)} u))$\\
  \hline
\end{tabular}
\endgroup
\end{center}\vspace{5pt}

In particular this reveals that, in the case of two distinct deciding factors, the diffusivity matrix is not enough to characterize heterogeneous diffusion. More precisely, the choice of two deciding factors allows us to recover the whole range of local diffusion equations that has been derived from kinetic equations in~\cite{KimSeo}. 

\subsection{Many or infinitely many deciding factors.}  The two above subsections  suggest that, as the number of deciding factors increase to account for more general heterogeneity, one may recover new local models through the singular limit procedure. In particular, it is worth pointing out that Wereide's law does not only appear as a limit of a single deciding 
factor \lq\lq close to departure'' (see subsection \ref{ss:one-DF}), but also as the singular limit of the focusing kernels associated with~\eqref{eq:strat}. We did not detail this computation since it proceeds similarly to what we have presented. In other words, the nonlocal Stratonovich model does converge to the Stratonovich type diffusion, as one may expect. We believe that generalizations of \eqref{eq:strat}, for which the {\it whole path decides}, deserve further investigations.

\bibliographystyle{siam}
\bibliography{biblio}

\end{document}